\documentclass[12pt]{article}
\usepackage[T1]{fontenc}
\usepackage{lmodern}
\usepackage{amsmath,amssymb,amsthm,mathtools}
\usepackage{mathrsfs}
\usepackage{microtype}
\usepackage[top=27mm, bottom=27mm, left=22mm, right=22mm]{geometry}
\usepackage{xcolor}
\usepackage{cite}
\usepackage[colorlinks=true,linkcolor=blue!50!black,citecolor=blue!50!black,urlcolor=blue!50!black]{hyperref}
\newtheorem{theorem}{Theorem}[section]
\newtheorem{corollary}[theorem]{Corollary}
\newtheorem{lemma}[theorem]{Lemma}
\newtheorem{question}[theorem]{Question}
\theoremstyle{definition}
\newtheorem{definition}[theorem]{Definition}

\newcommand{\cC}{\mathcal C}
\newcommand{\cT}{\mathcal T}
\newcommand{\cA}{\mathcal A}
\newcommand{\cL}{\mathcal L}
\newcommand{\cY}{\mathcal Y}
\newcommand{\Z}{\mathbb Z}
\newcommand{\N}{\mathbb N}
\newcommand{\E}{\mathbb E}
\newcommand{\Prob}{\mathbb P}
\DeclareMathOperator{\supp}{supp}
\DeclareMathOperator{\im}{Im}
\DeclareMathOperator{\Inj}{Inj}
\DeclareMathOperator{\Sym}{Sym}
\DeclareMathOperator{\id}{id}
\newcommand{\fall}[2]{(#1)_{#2}}
\newcommand{\one}{\mathbf 1}

\title{{\bf The existence of generalized Steiner systems}\thanks{The work of Shitao Li was supported by the National Natural Science Foundation of China under Grants 12526612 and 126011020. The work of Minjia Shi was supported by the National Natural Science Foundation of China under Grants 12471490 and W2611002.}}

\author{Xuan Wang\thanks{School of Mathematical Sciences, Anhui University, Hefei, Anhui 230601, China. E-mail: \textsf{wang\_xuan\_ah@163.com}.},~
Shitao Li\thanks{School of Internet, Anhui University, Hefei, Anhui 230039, China. E-mail: \textsf{lishitao0216@163.com}.},~
Minjia Shi\thanks{School of Mathematical Sciences, Anhui University, Hefei, Anhui 230601, China. E-mail: \textsf{smjwcl.good@163.com}.}}
\date{}

\begin{document}
\maketitle

\begin{abstract}
A generalized Steiner system $GS(t,w,n,q)$ is a $q$-ary constant-weight code of length $n$ and weight $w$ that has minimum distance at least $2(w-t)+1$ and attains the Johnson-type upper bound $(q-1)^t\binom nt/\binom wt$. For strengths $t=1$ and $t=2$, such systems were known to exist for all sufficiently large lengths satisfying the necessary divisibility conditions. For general fixed parameters, Liu and Shangguan (J. Combin.\ Theory Ser.~A 209 (2025), 105955) proved that the Johnson-type bound is attained asymptotically, and asked whether it is attained with equality for all sufficiently large lengths satisfying these conditions. We answer this question affirmatively. We prove that, for fixed integers $q\geq3$ and $w>t\geq1$, a $GS(t,w,n,q)$ exists for every sufficiently large $n$ such that $\binom{w-i}{t-i}$ divides $(q-1)^{t-i}\binom{n-i}{t-i}$ for all $0\leq i<t$. Moreover, the system can be chosen so that the supports of distinct codewords intersect in at most $t$ coordinates. Our proof assigns an independent random labelling to every set of $t+1$ coordinates, which turns large support intersections into conflicts in an exact covering problem, and then applies Keevash's decomposition theorem for labelled complexes.
\end{abstract}
\noindent\textbf{Keywords:} Constant-weight code, generalized Steiner system, Johnson bound, labelled complex, design existence.

\noindent\textbf{2020 MSC:} Primary 05B05, 94B25; Secondary 05B30, 05C70, 05D40.

\section{Introduction}\label{sec:intro}

A classical problem in combinatorial design theory, going back to Steiner \cite{Steiner1853}, is to determine when a Steiner system exists. For integers $1\leq t<w<n$, a \emph{Steiner system} with these parameters is a family of $w$-subsets of an $n$-set, called \emph{blocks}, such that every $t$-subset is contained in exactly one block. Counting the blocks that contain a fixed $i$-subset shows that such a system can exist only if $\binom{w-i}{t-i}$ divides $\binom{n-i}{t-i}$ for every $0\leq i<t$, and the existence conjecture asserts that these divisibility conditions suffice for all sufficiently large $n$. In 1963, Erd\H{o}s and Hanani \cite{ErdosHanani1963} conjectured an approximate version, namely that for fixed $t<w$ one can choose $(1-o(1))\binom nt/\binom wt$ blocks such that every $t$-subset is contained in at most one of them. R\"odl \cite{Rodl1985} proved this conjecture in 1985 by introducing the nibble method. Strengthening R\"odl's result, Keevash \cite{Keevash2014} proved the existence conjecture itself in 2014, and Glock, K\"uhn, Lo and Osthus \cite{GlockKuhnLoOsthus2023} later gave a different proof via iterative absorption, which also applies to decompositions into copies of an arbitrary uniform hypergraph.

Steiner systems also arise naturally in coding theory. Let $[n]=\{1,\ldots,n\}$, and identify a $w$-subset of $[n]$ with its characteristic vector, which is a binary word of length $n$ and weight $w$. Two $w$-subsets sharing exactly $s$ points correspond to words at Hamming distance $2(w-s)$. Hence, a family of $w$-subsets in which any two members share at most $t-1$ points is the same as a binary code of constant weight $w$ and minimum distance at least $2(w-t)+1$. Since every $t$-subset is contained in at most one member of such a family, the code has at most $\binom nt/\binom wt$ words, with equality if and only if the family is a Steiner system. In this language, R\"odl's theorem determines the maximum size of such a code asymptotically, and Keevash's theorem shows that the bound $\binom nt/\binom wt$ is attained for all sufficiently large $n$ satisfying the divisibility conditions.

It is natural to ask what happens over larger alphabets. Throughout this paper, $q$ denotes the size of the alphabet $\mathbb{Z}_q=\{0,1,\ldots,q-1\}$, in which $0$ is distinguished, and we put $m=q-1$. For $x\in\mathbb{Z}_q^n$, the \emph{support} of $x$ is $\supp(x)=\{i\in[n]:x_i\ne0\}$, and the \emph{weight} of $x$ is $|\supp(x)|$. The Hamming distance $d_H(x,y)$ is the number of coordinates at which $x$ and $y$ differ. Let $A_q(n,d,w)$ be the maximum size of a code in $\mathbb{Z}_q^n$ in which every word has weight $w$ and any two distinct words have Hamming distance at least $d$. For $q\geq3$, a word of weight $w$ is no longer determined by its support, and two words may share coordinates at which their nonzero symbols differ. The role of $t$-subsets is now played by words of weight $t$, and inclusion of sets is replaced by the following notion. We say that a word $x$ \emph{covers} a word $z$ if $\supp(z)\subseteq\supp(x)$ and $x_i=z_i$ for every $i\in\supp(z)$. For an odd distance $d$ with $1\leq d\leq 2w+1$, put $t=(2w-d+1)/2$, so that $d=2(w-t)+1$ and $0\leq t\leq w$. The Johnson-type bound \cite[Theorem~1.2]{LiuShangguan2025} states that
\begin{equation}\label{eq:johnson}
A_q(n,2w-2t+1,w) \leq \frac{(q-1)^t\binom nt}{\binom wt}.
\end{equation}
For $q=2$ and $1\leq t<w$, the bound \eqref{eq:johnson} reduces to the bound $\binom nt/\binom wt$ discussed above. The same theorem gives $A_q(n,2(w-t)+2,w)\leq(q-1)^{t-1}\binom nt/\binom wt$ for $1\leq t\leq w$ in the even case, which is not needed here.

\subsection{Generalized Steiner systems}\label{subsec:GS}

In his study of optimal constant-weight codes over arbitrary alphabets, Etzion \cite{Etzion1997} introduced the following generalization of Steiner systems, which we state in our notation.

\begin{definition}[Generalized Steiner system]\label{def:GS}
Let $q\geq3$ and $1\leq t<w\leq n$. A code $\cC\subseteq\mathbb{Z}_q^n$ whose codewords all have weight $w$ and whose distinct codewords have Hamming distance at least $2(w-t)+1$ is a \emph{generalized Steiner system} $GS(t,w,n,q)$ if every $q$-ary word of weight $t$ is covered by exactly one codeword of $\cC$.
\end{definition}

By Lemma~\ref{lem:counting} below, a code of weight $w$ and minimum distance at least $2(w-t)+1$ is a $GS(t,w,n,q)$ if and only if it attains the upper bound in \eqref{eq:johnson}. The same lemma shows, by counting the codewords that cover a fixed word of weight $i$, that a $GS(t,w,n,q)$ can exist only if
\begin{equation}\label{eq:admissible}
\binom{w-i}{t-i}\ \bigm|\ (q-1)^{t-i}\binom{n-i}{t-i}
\ \text{ for all } 0\leq i<t.
\end{equation}
For fixed $q,w,t$, a length $n$ satisfying \eqref{eq:admissible} is called \emph{admissible}.

In \cite{Etzion1997}, the definition is formulated in terms of H-designs. The points of such a design are partitioned into $n$ groups of size $g$, its blocks are $w$-sets meeting every group in at most one point, and every $t$-set with this property is contained in exactly one block. If the points of the $i$th group are identified with the nonzero symbols at coordinate $i$, then each block becomes a word of weight $w$ over an alphabet of size $g+1$. A generalized Steiner system is an H-design whose associated code has minimum distance at least $2(w-t)+1$, and its notation uses the group size $g$ as the last parameter. This is Definition~\ref{def:GS} with $q=g+1$. Thus, a system whose last parameter is $g$ in \cite{Etzion1997} and in much of the design literature has last parameter $q=g+1$ in this paper. We also note that the term $q$-Steiner system is used in finite geometry for designs whose blocks are subspaces, and such designs are not considered here.

Generalized Steiner systems are also extremal for the code--anticode bound. The words of weight $w$ covering a fixed word of weight $t$ form an anticode of size $(q-1)^{w-t}\binom{n-t}{w-t}$ and diameter at most $2(w-t)$, and the product of this size with the right-hand side of \eqref{eq:johnson} equals the number of words of weight $w$. For $n\geq 2w-t$, the anticode has diameter exactly $2(w-t)$, so every $GS(t,w,n,q)$ is $2(w-t)$-diameter perfect and attains the code--anticode bound \cite[Lemma~11 and Corollary~5]{Etzion2022}, see also \cite[Section~9.3]{Etzion2022Book}. Moreover, for $t>1$, fixing a coordinate and a nonzero symbol, retaining the codewords with that symbol at the chosen coordinate, and deleting that coordinate gives a $GS(t-1,w-1,n-1,q)$, see \cite{Etzion2022Book}. To see this, note that the retained words agree at the deleted coordinate, so deletion preserves their pairwise distances. In addition, a word of weight $t-1$ on the remaining coordinates, extended by the chosen symbol at the deleted coordinate, becomes a word of weight $t$, and the unique codeword covering it is retained.

Etzion's constructions in \cite{Etzion1997} connect such codes with block designs and include families obtained from ordinary Steiner systems on the supports. Since then, generalized Steiner systems of small strength and weight have been studied extensively. Generalized Steiner triple systems, written $GS(2,3,n,g)$ in the group-size notation, correspond to $(g+1)$-ary constant-weight codes of length $n$, weight $3$ and minimum distance at least $3$ that attain the Johnson-type bound. Phelps and Yin \cite{PhelpsYin1997} gave constructions for $g\equiv 3\pmod 6$. Blake-Wilson and Phelps \cite{BlakeWilsonPhelps1999} proved that, for every fixed group size $g$, the necessary conditions for the existence of a $GS(2,3,n,g)$ are sufficient for all sufficiently large $n$. A complementary line of work determined the full existence spectrum for particular group sizes, including $g=4$ \cite{PhelpsYin1999}, $g=5,6$ \cite{ChenGeZhu1999,ChenGeZhu2000}, and $g=7,8$ \cite{WuGeZhu2000}. Ge \cite{Ge2000,Ge2002ZeroThree,Ge2002OneFive} developed finite-range reductions and established the sufficiency of the necessary conditions for $g=11,12,13,15$. For a survey of strength-two systems and related constructions via maximum distance holey packings, see Yin \cite{Yin2002}. Group divisible designs have also been used more broadly to construct optimal nonbinary constant-weight codes, including the quaternary families of Wu and Fan \cite{WuFan2009}. For strength three, Cao, Ji and Zhu \cite{CaoJiZhu2007} developed constructions for ternary systems of block size four, denoted $GS(3,4,v,2)$ in the group-size notation. A related line of work asks for many disjoint exact systems. Tan and Zhou \cite{TanZhou2026} studied large sets of generalized Steiner systems as tilings of the entire constant-weight space, and their constructions give infinite families of such tilings for selected parameters, including nonbinary weight-three cases.

For arbitrary fixed alphabet size and weight, two regimes were known in which \eqref{eq:johnson} is attained for all sufficiently large admissible lengths. When $t=1$, the target distance is $2w-1$, and \eqref{eq:admissible} reads $w\mid(q-1)n$. Etzion \cite[Problem~7]{Etzion1997} asked whether, for fixed $q$ and $w$, a $GS(1,w,n,q)$ exists for every sufficiently large $n$ satisfying this condition. Chee, Dau, Ling and Ling \cite[Main Theorem~2 and Theorem~6.1]{CheeDauLingLing2010} answered this question affirmatively by constructing quasicyclic codes of size $(q-1)n/w$ using generalized difference triangle sets. When $t=2$, the target distance is $2w-3$, and \eqref{eq:admissible} consists of the two conditions
\[
w(w-1)\mid(q-1)^2n(n-1)
\ \text{ and }\
w-1\mid(q-1)(n-1).
\]
Chee et al.\ \cite[Proposition~4.4]{CheeEtAl2019} proved, via superpure decompositions of edge-colored digraphs, that \eqref{eq:johnson} is attained whenever $n$ is sufficiently large and satisfies these two conditions.

For general fixed $q$, $w$ and $t$, Liu and Shangguan \cite[Theorem~1.4]{LiuShangguan2025} recently established the asymptotic formula
\[
A_q(n,2w-2t+1,w) =(1-o(1))\frac{(q-1)^t\binom nt}{\binom wt},
\]
where $o(1)\to0$ as $n\to\infty$. This is a nonbinary counterpart of R\"odl's theorem, and it implies the existence of near-optimal generalized Steiner systems. Their proof for constant-weight codes uses a near-perfect matching in an auxiliary hypergraph. For constant-composition codes, they also used a conflict-free matching approach, forbidding pairs whose supports intersect in more than $t$ coordinates \cite[Lemma~4.1]{LiuShangguan2025}. Asymptotic tightness, however, allows a vanishing proportion of the words of weight $t$ to remain uncovered, whereas a generalized Steiner system must cover each of them exactly once. Motivated by Keevash's theorem in the binary case, Liu and Shangguan raised the following question.

\begin{question}[Liu--Shangguan {\cite[Question~5.3]{LiuShangguan2025}}]\label{q:LS}
Let $q\geq3$, and let $d$ be odd. Can the Johnson-type bound for $A_q(n,d,w)$ be attained with equality for all sufficiently large $n$ satisfying the necessary divisibility conditions?
\end{question}

With $d=2(w-t)+1$, the necessary divisibility conditions are those in \eqref{eq:admissible}, so Question~\ref{q:LS} asks whether a $GS(t,w,n,q)$ exists for every sufficiently large admissible $n$. Further background and open problems on generalized Steiner systems can be found in Etzion's book \cite{Etzion2022Book}.

\subsection{Main results}\label{subsec:main}

In this paper, we answer Question~\ref{q:LS} affirmatively. Our main result is the following.

\begin{theorem}\label{thm:main}
Let $q\geq 3$ and $w>t\geq 1$ be fixed integers. Then there is an integer $n_0=n_0(q,w,t)$ such that a $GS(t,w,n,q)$ exists for every $n\geq n_0$ satisfying \eqref{eq:admissible}.
\end{theorem}

Combining Lemma~\ref{lem:counting} with Theorem~\ref{thm:main}, we obtain the following characterization.

\begin{corollary}\label{cor:characterization}
Let $q\geq 3$ and $w>t\geq 1$ be fixed integers. For all sufficiently large $n$, a $GS(t,w,n,q)$ exists if and only if $n$ satisfies \eqref{eq:admissible}. In particular, $A_q(n,2w-2t+1,w)=(q-1)^t\binom nt/\binom wt$ for all sufficiently large admissible $n$.
\end{corollary}

The two endpoint cases of Question~\ref{q:LS} are immediate. For $t=w$, the distance is $1$, and one takes all words of weight $w$. For $t=0$, the required distance is $2w+1$, and an optimal code consists of a single word, since any two words of weight $w$ are at distance at most $2w$. Together with Theorem~\ref{thm:main}, these two cases answer Question~\ref{q:LS} affirmatively for every odd distance $d\leq 2w+1$. For $t=1$ and $t=2$, Theorem~\ref{thm:main} includes the existence results of \cite{CheeDauLingLing2010} and \cite{CheeEtAl2019} mentioned above. Just as the result of Liu and Shangguan is a nonbinary counterpart of R\"odl's theorem, Theorem~\ref{thm:main} can be viewed as a nonbinary counterpart of Keevash's theorem on the existence of Steiner systems.

\subsection{Our approach}\label{subsec:approach}

We work with labelled sets, which record the coordinates and the symbols of a word separately. For $0\leq j\leq n$, let
\[
\cT_j(n,m)=\left\{(S,a):S\in\binom{[n]}{j},\ a:S\longrightarrow[m]\right\}.
\]
An element $(S,a)\in\cT_j(n,m)$ is a \emph{labelled $j$-set}, and we write $(S,a)\subseteq(W,b)$ if $S\subseteq W$ and $b|_S=a$. A word $x$ of weight $j$ corresponds to the labelled $j$-set $(\supp(x),x|_{\supp(x)})$. Conversely, a labelled $j$-set determines a word of weight $j$ by placing zeros outside its support. Under this correspondence, $x$ covers $z$ if and only if the labelled set of $z$ is contained in that of $x$. Thus, the covering condition in Definition~\ref{def:GS} says that every labelled $t$-set in $\cT_t(n,m)$ is contained in exactly one of the labelled $w$-sets corresponding to the codewords.

A labelled set $(S,a)$ may also be viewed as the set $\{(i,a(i)):i\in S\}$ of coordinate--symbol pairs, which is a transversal subset of $[n]\times[m]$, see \eqref{eq:label-transversal}. In this representation, the covering condition asks for an exact cover of the transversal $t$-subsets of $[n]\times[m]$ by transversal $w$-subsets, which is the H-design condition in Etzion's formulation. The main difficulty lies in the distance condition. By \eqref{eq:distance}, the distance between two words of weight $w$ is determined by the number of coordinate--symbol pairs they share and by the size of the intersection of their supports. Exact coverage ensures that two distinct codewords share at most $t-1$ coordinate--symbol pairs, but it does not bound the number of common coordinates at which their symbols differ. We therefore also require that the supports of distinct codewords intersect in at most $t$ coordinates, and then \eqref{eq:distance} gives the required minimum distance, see Lemma~\ref{lem:distance-transfer}. Theorem~\ref{thm:main} will be deduced from the following labelled covering result.

\begin{theorem}\label{thm:labelled}
Let $q\geq3$ and $w>t\geq1$ be fixed integers, and put $m=q-1$. For every sufficiently large $n$ satisfying \eqref{eq:admissible}, there exists a family $\mathcal B\subseteq\cT_w(n,m)$ with the following properties.
\begin{enumerate}
\item[(i)] Every member of $\cT_t(n,m)$ is contained in exactly one member of $\mathcal B$.
\item[(ii)] If $(W,a)$ and $(W',a')$ are distinct members of $\mathcal B$, then $|W\cap W'|\leq t$.
\end{enumerate}
\end{theorem}

A family as in Theorem~\ref{thm:labelled} also determines an ordinary design. A coordinate $t$-set $T$ has $m^t$ labellings, each of which is contained in exactly one block of the family, and a block whose support contains $T$ realizes exactly one of them. Hence, the supports of the blocks form a $t$-$(n,w,m^t)$ design. Since $w>t$, property~(ii) ensures that distinct blocks have distinct supports, so this design is simple, and any two of its blocks intersect in at most $t$ points. Property~(ii) also gives a lower bound on the length. The $m^t$ labellings of $T$ are covered by distinct blocks whose supports pairwise intersect exactly in $T$, so the parts of these supports outside $T$ are pairwise disjoint. Thus, every family as in Theorem~\ref{thm:labelled} has $n\geq t+m^t(w-t)$. This lower bound concerns only exact covers with the support property, and it does not affect Theorem~\ref{thm:labelled}, in which $q,w,t$ are fixed and $n$ is large.

To prove Theorem~\ref{thm:labelled}, we prescribe a random labelling $\chi_S:S\to[m]$ of every coordinate $(t+1)$-set $S$, independently for different sets, and allow only blocks whose restriction to each $(t+1)$-subset $S$ of their support is $\chi_S$. If the supports of two allowed blocks share a set $S$ of $t+1$ coordinates, then both blocks carry the labelling $\chi_S$ on $S$, and hence they contain a common labelled $t$-set. Consequently, in an exact labelled cover by allowed blocks, distinct supports meet in at most $t$ coordinates. Finding such a cover is a decomposition problem for a labelled complex. Keevash \cite{Keevash2018} generalized the existence of designs to exact decompositions of vectors indexed by the labelled faces of simplicial complexes, and we apply his theorem to the random complex defined by these restrictions. Its three hypotheses are verified as follows. Extendability follows from a concentration estimate for one-vertex extensions. Integer feasibility is transferred from the complete complex, which has the same faces of size at most $t$ as the random complex. Fractional regularity follows from explicit weights that balance the contributions of blocks with different compositions.

The paper is organized as follows. In Section~\ref{sec:prelim}, we translate the distance and covering conditions of a generalized Steiner system into conditions on labelled sets, and we recall Keevash's decomposition theorem for labelled complexes in the form needed here. In Section~\ref{sec:proof}, we prove Theorem~\ref{thm:labelled}, and hence Theorem~\ref{thm:main}, assuming three lemmas on the random complex described above. A proof overview is given at the beginning of that section. The three lemmas are proved in Section~\ref{sec:lemmas}. We conclude the paper with some open questions in Section~\ref{sec:conclusion}.

\section{Preliminaries}\label{sec:prelim}

\subsection{Constant-weight codes and labelled coverings}\label{subsec:codes}

In this subsection, we translate the distance and covering conditions of a generalized Steiner system into conditions on labelled sets. Throughout, $m=q-1$, and words of weight $w$ are identified with labelled $w$-sets as in Section~\ref{subsec:approach}.

Let $X=[n]\times[m]$, and call a subset of $X$ \emph{transversal} if its elements have distinct first coordinates. The map
\begin{equation}\label{eq:label-transversal}
u=(S,a)\longmapsto \widehat u=\{(i,a(i)):i\in S\}
\end{equation}
is a bijection from $\cT_j(n,m)$ to the family of transversal $j$-subsets of $X$. For labelled sets $u$ and $W$, we have $u\subseteq W$ if and only if $\widehat u\subseteq\widehat W$. For two words $x$ and $y$ of weight $w$, put
\[
s(x,y)=|\supp(x)\cap\supp(y)|
\ \text{ and }\
a(x,y)=|\{i:x_i=y_i\ne0\}|,
\]
so that $a(x,y)=|\widehat x\cap\widehat y|$, where $\widehat x$ denotes the transversal set of the labelled set of $x$. Exactly $2(w-s(x,y))$ coordinates lie in only one of the two supports, and exactly $s(x,y)-a(x,y)$ common coordinates carry different symbols. Hence, the Hamming distance between $x$ and $y$ is
\begin{equation}\label{eq:distance}
d_H(x,y)=2w-s(x,y)-a(x,y).
\end{equation}
The next lemma derives the bound \eqref{eq:johnson} and the divisibility conditions \eqref{eq:admissible} by counting labelled sets.

\begin{lemma}\label{lem:counting}
Let $\cC$ be a code of length $n$ and weight $w$ with minimum distance at least $2w-2t+1$. Then each labelled $t$-set is contained in at most one word of $\cC$. Consequently, \eqref{eq:johnson} holds, and $|\cC|$ equals the right-hand side of \eqref{eq:johnson} if and only if $\cC$ is a $GS(t,w,n,q)$. In this case, for every labelled $i$-set $u$ with $0\leq i<t$, the number of words of $\cC$ containing $u$ is
\begin{equation}\label{eq:lambda-i}
\lambda_i=\frac{m^{t-i}\binom{n-i}{t-i}}{\binom{w-i}{t-i}}.
\end{equation}
In particular, if a $GS(t,w,n,q)$ exists, then $n$ satisfies \eqref{eq:admissible}.
\end{lemma}

\begin{proof}
Suppose that two distinct words $x,y\in\cC$ contain the same labelled $t$-set. Then $s(x,y)\geq t$ and $a(x,y)\geq t$, and \eqref{eq:distance} gives $d_H(x,y)\leq 2w-2t$, a contradiction. Each word contains $\binom wt$ labelled $t$-sets, while there are $m^t\binom nt$ labelled $t$-sets in total. This proves \eqref{eq:johnson}, and equality holds precisely when every labelled $t$-set is contained in exactly one word, that is, when $\cC$ is a $GS(t,w,n,q)$.

Now assume equality, and fix a labelled $i$-set $u$ with $0\leq i<t$. There are $m^{t-i}\binom{n-i}{t-i}$ labelled $t$-sets containing $u$. Each of them lies in exactly one word, and this word contains $u$. Conversely, each word containing $u$ contains exactly $\binom{w-i}{t-i}$ of these labelled $t$-sets. Counting the pairs formed by such a word and such a labelled $t$-set in two ways gives \eqref{eq:lambda-i}. In particular, all $\lambda_i$ are integers, which is \eqref{eq:admissible}.
\end{proof}

The next lemma shows that a labelled covering with small support intersections defines a generalized Steiner system.

\begin{lemma}\label{lem:distance-transfer}
Let $m\geq2$ and $1\leq t<w\leq n$, and let $\mathcal B\subseteq\cT_w(n,m)$. Suppose that every labelled $t$-set is contained in at most one member of $\mathcal B$, and that distinct members have supports intersecting in at most $t$ coordinates. Then the associated $(m+1)$-ary constant-weight code has minimum distance at least $2w-2t+1$. If every labelled $t$-set is covered exactly once, then this code is a $GS(t,w,n,m+1)$ of size $m^t\binom nt/\binom wt$, and its minimum distance is exactly $2w-2t+1$.
\end{lemma}

\begin{proof}
For distinct codewords $x,y$, the uniqueness of covering implies $a(x,y)\leq t-1$, and the support condition implies $s(x,y)\leq t$. By \eqref{eq:distance},
\[
d_H(x,y)=2w-s(x,y)-a(x,y)\geq2w-2t+1.
\]
If the coverage is exact, the code therefore satisfies the definition of a generalized Steiner system. Counting incidences between blocks and labelled $t$-sets gives $|\mathcal B|\binom wt=m^t\binom nt$, which is the claimed size.

To show that the distance bound is attained, fix a coordinate $t$-set $T$ and two maps $a,a':T\to[m]$ that differ at exactly one coordinate. Such maps exist because $m\geq2$ and $t\geq1$. Let $x$ and $y$ be the codewords covering $(T,a)$ and $(T,a')$, respectively. They are distinct, because no word induces both labellings on $T$. They share the $t$ coordinates of $T$ and agree on $t-1$ of them, so $d_H(x,y)\leq2w-t-(t-1)$. Together with the lower bound, this proves equality.
\end{proof}

\subsection{Labelled complexes and Keevash's decomposition theorem}\label{subsec:keevash}

In this subsection, we recall the framework of Keevash \cite{Keevash2018} in the form needed here and derive the decomposition statement used in Section~\ref{sec:proof}. We first introduce complexes and their symmetries, then the extendability condition, and finally the templates and molecules that encode the codewords. These notions lead to Keevash's theorem and to the special case used in our proof, Lemma~\ref{lem:keevash-special}. All definitions can be found in \cite[Sections~2 and~3]{Keevash2018}, and all numbered results of \cite{Keevash2018} refer to its first arXiv version, arXiv:1802.05900v1.

For the remainder of the paper, fix $q\geq3$ and $w>t\geq1$, and put $m=q-1$ and $L=mw$. Let $R=[m]\times[w]$ be the set of \emph{roles}, partitioned as $R=P_1\sqcup\cdots\sqcup P_m$ with $P_a=\{a\}\times[w]$. Write $\kappa(u)=a$ for $u\in P_a$, and call $\kappa(u)$ the \emph{color} of $u$. The second index of a role $(a,j)$ distinguishes copies of the symbol $a$ and is not a coordinate of the code. Let $\Sigma=\Sym(P_1)\times\cdots\times\Sym(P_m)$ be the group of color-preserving permutations of $R$. Every color-preserving bijection between two subsets of $R$ extends to an element of $\Sigma$, since it can be extended separately within each class $P_a$.

The basic objects of \cite{Keevash2018} are complexes of injections defined on subsets of $R$. For a set $V$ and $B\subseteq R$, let $\Inj(B,V)$ be the set of all injections from $B$ to $V$. An $R$-\emph{complex} $\Phi$ on $V$ is a family $\Phi=(\Phi_B:B\subseteq R)$ with $\Phi_B\subseteq\Inj(B,V)$ that is closed under restriction, that is, $\phi|_C\in\Phi_C$ whenever $\phi\in\Phi_B$ and $C\subseteq B$. Each $\phi\in\Phi_B$ is a \emph{face} of $\Phi$ with domain $B$ and size $|B|$. The vertex set $V(\Phi)$ consists of all $\phi(u)$ with $u\in R$ and $\phi\in\Phi_{\{u\}}$. The \emph{complete} $R$-complex $\Phi^0$ on $[n]$ is given by $\Phi^0_B=\Inj(B,[n])$ for every $B\subseteq R$. An $R$-complex $\Psi$ with $\Psi_B\subseteq\Phi_B$ for every $B\subseteq R$ is a \emph{subcomplex} of $\Phi$. For faces $\psi'$ and $\psi$, we write $\psi'\subseteq\psi$ if $\psi'$ is a restriction of $\psi$. For a finite set $I$, we write $\Z^I$ for the set of functions $I\to\Z$, and $\one_O\in\Z^I$ for the indicator vector of $O\subseteq I$.

The group $\Sigma$ acts by relabelling roles. For $\psi\in\Phi_B$ and $\sigma\in\Sigma$, define
\[
\psi\sigma=\psi\circ\bigl(\sigma|_{\sigma^{-1}(B)}\bigr):
\sigma^{-1}(B)\longrightarrow V.
\]
This operation changes the domain but preserves the image, and since $\sigma$ preserves colors, the role occupying a given image vertex keeps its color. We say that $\Phi$ is $\Sigma$-\emph{adapted} if
\begin{equation}\label{eq:adapted}
\psi\sigma\in\Phi_{\sigma^{-1}(B)}
\ \text{ for all } B\subseteq R,\ \psi\in\Phi_B \text{ and } \sigma\in\Sigma.
\end{equation}
In this case, the \emph{orbit} of a face $\psi$ is $\psi\Sigma=\{\psi\sigma:\sigma\in\Sigma\}$. For $0\leq j\leq L$, we write $\Phi_j$ for the set of faces of size $j$, $\Phi_j/\Sigma$ for the set of their orbits, and $\Phi/\Sigma$ for the set of all face orbits.

A role $t$-face $\psi\in\Phi_t$ of a complex on $[n]$ induces a labelled coordinate $t$-set. Define $a_\psi:\im\psi\to[m]$ by $a_\psi(i)=\kappa(\psi^{-1}(i))$, and put
\[
\ell(\psi)=(\im\psi,a_\psi)\in\cT_t(n,m).
\]
Thus, the coordinate $\psi(u)$ receives the symbol $\kappa(u)$ of its role $u$, while $\ell(\psi)$ forgets the name of that role. Every $R$-complex on $[n]$ considered in this paper satisfies the \emph{complete skeleton condition}
\begin{equation}\label{eq:full-skeleton}
\Phi_B=\Phi_B^0=\Inj(B,[n])
\ \text{ for all } B\subseteq R \text{ with } |B|\leq t.
\end{equation}

\begin{lemma}\label{lem:orbit-identification}
Let $\Phi$ be a $\Sigma$-adapted $R$-complex on $[n]$ satisfying \eqref{eq:full-skeleton}. For any $\psi,\psi'\in\Phi_t$, we have $\ell(\psi)=\ell(\psi')$ if and only if $\psi'\in\psi\Sigma$. Consequently, $\ell$ induces a bijection
\begin{equation}\label{eq:orbit-identification}
\Phi_t/\Sigma\longrightarrow\cT_t(n,m),
\ \ \psi\Sigma\longmapsto\ell(\psi).
\end{equation}
\end{lemma}

\begin{proof}
Suppose first that $\psi'=\psi\sigma$ for some $\sigma\in\Sigma$. Then $\psi'$ has the same image as $\psi$. At the coordinate $i=\psi'(u)=\psi(\sigma(u))$, the two faces induce the symbols $\kappa(u)$ and $\kappa(\sigma(u))$, which are equal because $\sigma$ preserves colors. Hence, $\ell(\psi')=\ell(\psi)$.

Conversely, suppose that $\ell(\psi')=\ell(\psi)$, where $\psi\in\Phi_B$ and $\psi'\in\Phi_{B'}$. Since the two faces have the same image, $\beta=\psi^{-1}\circ\psi':B'\to B$ is a well-defined bijection. For each $u\in B'$, the equality of the induced labels at $\psi'(u)$ gives $\kappa(\beta(u))=a_\psi(\psi'(u))=a_{\psi'}(\psi'(u))=\kappa(u)$. Thus, $\beta$ preserves colors and extends to some $\sigma\in\Sigma$ with $\sigma|_{B'}=\beta$. Since $\sigma(B')=B$, the face $\psi\sigma$ has domain $B'$ and satisfies $(\psi\sigma)(u)=\psi(\beta(u))=\psi'(u)$ for $u\in B'$. Therefore, $\psi'=\psi\sigma$.

It remains to show that every $(S,a)\in\cT_t(n,m)$ is of the form $\ell(\psi)$. For each $k\in[m]$, put $c_k=|a^{-1}(k)|$. Since $c_k\leq t<w=|P_k|$, we can choose $B_k\subseteq P_k$ of size $c_k$ and a bijection $\psi_k:B_k\to a^{-1}(k)$. The sets $B_k$ are pairwise disjoint and the sets $a^{-1}(k)$ partition $S$, so the maps $\psi_k$ combine into a bijection $\psi:B\to S$, where $B=B_1\sqcup\cdots\sqcup B_m$ has size $t$. By \eqref{eq:full-skeleton}, we have $\psi\in\Inj(B,[n])=\Phi_B$. Moreover, $a(\psi(u))=\kappa(u)$ for every $u\in B$ and $\im\psi=S$, so $\ell(\psi)=(S,a)$.
\end{proof}

We next recall the extendability condition, which requires every partial embedding of a bounded complex to have many completions. Let $H$ and $\Phi$ be $R$-complexes. A $\Phi$-\emph{embedding} of $H$ is an injective map $g:V(H)\to V(\Phi)$ such that $g\circ\theta\in\Phi_B$ for every $B\subseteq R$ and every $\theta\in H_B$. For $h\in\N$, let $R(h)$ be the $R$-complex on $R\times[h]$ whose faces with domain $B$ are the \emph{partite} maps $\theta:B\to R\times[h]$, namely those with $\theta(u)=(u,j_u)$ for every $u\in B$, where $j_u\in[h]$. The vertex $(u,j)$ is given the color $\kappa_h(u,j)=\kappa(u)$. For a subcomplex $H\subseteq R(h)$ and $F\subseteq V(H)$, let $H[F]$ consist of the faces of $H$ whose images lie in $F$. Every vertex of $H$ lies in a singleton face, so $V(H[F])=F$. A \emph{rooted extension} is a triple $(H,F,f)$, where $f:F\to V(\Phi)$ is a $\Phi$-embedding of $H[F]$, and a \emph{completion} of it is a $\Phi$-embedding $g$ of $H$ with $g|_F=f$.

\begin{definition}[see {\cite[Definition~2.9]{Keevash2018}}]\label{def:extendability}
Let $0<\omega\leq1$ and $h\in\N$. An $R$-complex $\Phi$ is $(\omega,h)$-\emph{extendable} if, for every subcomplex $H\subseteq R(h)$, every $F\subseteq V(H)$ and every $\Phi$-embedding $f$ of $H[F]$, the number of completions of $(H,F,f)$ is at least $\omega\,|V(\Phi)|^{|V(H)\setminus F|}$. Since $R(h')\subseteq R(h)$ for $h'\leq h$, an $(\omega,h)$-extendable complex is also $(\omega',h')$-extendable whenever $0<\omega'\leq\omega$ and $1\leq h'\leq h$.
\end{definition}

We now describe the templates, whose embeddings produce the candidate codewords, and the vectors that these embeddings contribute to a decomposition. Let $\Sigma^{\leq}$ be the $R$-complex of all restrictions of permutations in $\Sigma$. Its faces with domain $B$ are exactly the color-preserving injections $B\to R$, since every such injection extends to an element of $\Sigma$. For each $J\in\binom Rw$, let $A_J$ be a separately indexed copy of $\Sigma^{\leq}$, called a \emph{template}, whose vertex set we identify with $R$. The roles in $J$ are \emph{active}, and those in $R\setminus J$ are \emph{auxiliary}. Formally, a template face is an indexed pair $(J,\theta)$, so that templates with different active sets are distinct although their underlying complexes coincide. Let $\cA=\{A_J:J\in\binom Rw\}$. A $\Phi$-embedding of $A_J$ is an injective map $\phi:R\to[n]$ with $\phi\circ\theta\in\Phi_B$ for every $B\subseteq R$ and every $\theta\in(A_J)_B$. Write $A_J(\Phi)$ for the set of these maps, and
\[
\cA(\Phi)=\left\{(J,\phi):J\in\binom Rw,\ \phi\in A_J(\Phi)\right\}
\]
for the set of \emph{template embeddings}. For $B\subseteq R$, let $\cA_B$ be the disjoint union of the sets $(A_J)_B$ over all $J$. These notions specialize \cite[Definitions~2.19 and~2.21]{Keevash2018}.

\begin{lemma}[see {\cite[Lemma~2.20]{Keevash2018}}]\label{lem:embedding-identification}
Let $\Phi$ be a $\Sigma$-adapted $R$-complex. Then $A_J(\Phi)=\Phi_R$ for every $J\in\binom Rw$.
\end{lemma}

\begin{proof}
If $\phi\in A_J(\Phi)$, then $\phi=\phi\circ\id_R\in\Phi_R$, since $\id_R$ is a face of $A_J$. Conversely, let $\phi\in\Phi_R$ and $\theta\in(A_J)_B$. Then $\theta$ extends to some $\sigma\in\Sigma$ with $\sigma|_B=\theta$. Since $\Phi$ is closed under restriction, we have $\phi|_{\sigma(B)}\in\Phi_{\sigma(B)}$, and since $\Phi$ is $\Sigma$-adapted, $\phi\circ\theta=(\phi|_{\sigma(B)})\circ(\sigma|_B)\in\Phi_B$. Hence, $\phi\in A_J(\Phi)$.
\end{proof}

Although $A_J(\Phi)$ does not depend on $J$, we keep the pair $(J,\phi)$, since the same full map produces different codewords for different active sets. The \emph{active word} of a template embedding $(J,\phi)$ is the labelled $w$-set $x_{J,\phi}=(\phi(J),a_{J,\phi})$ with $a_{J,\phi}(\phi(u))=\kappa(u)$ for $u\in J$. Viewed as a word of length $n$, it is
\begin{equation}\label{eq:active-word}
(x_{J,\phi})_i=
\begin{cases}
\kappa(\phi^{-1}(i)),&i\in\phi(J),\\
0,&i\notin\phi(J).
\end{cases}
\end{equation}
The auxiliary roles contribute no nonzero coordinates to this word, but they must satisfy the embedding conditions. For each $J$ and each $t$-face $\theta:B\to R$ of $A_J$, define
\begin{equation}\label{eq:gamma}
\gamma_{J,\theta}=\one_{\im\theta\subseteq J}.
\end{equation}
Here it is the image $\theta(B)$, rather than the domain $B$, that is required to lie in $J$. The \emph{molecule} of $(J,\phi)$ is the vector $v_{J,\phi}\in\Z^{\Phi_t}$ obtained by transporting these coefficients to the host, that is, for $\psi\in\Phi_B$ with $|B|=t$,
\begin{equation}\label{eq:molecule}
v_{J,\phi}(\psi)=
\begin{cases}
\gamma_{J,\theta},&\text{if }\psi=\phi\circ\theta\text{ for some }
\theta\in(A_J)_B,\\
0,&\text{if no such face exists}.
\end{cases}
\end{equation}
Since $\phi$ is injective, any such $\theta$ equals $\phi^{-1}\circ\psi$, so the molecule is well defined. In particular, $v_{J,\phi}(\psi)=1$ if and only if $\psi=\phi\circ\theta$ for a color-preserving injection $\theta:B\to J$. We also write $v_e=v_{J,\phi}$ for $e=(J,\phi)$.

\begin{lemma}\label{lem:molecule-word}
Let $J\in\binom Rw$, $\phi\in\Phi_R$, and $\psi\in\Phi_B$ with $|B|=t$. Then
\begin{equation}\label{eq:molecule-word}
v_{J,\phi}(\psi)=1
\ \Longleftrightarrow\
\ell(\psi)\subseteq x_{J,\phi}.
\end{equation}
\end{lemma}

\begin{proof}
Suppose first that $v_{J,\phi}(\psi)=1$. Then $\psi=\phi\circ\theta$ for a color-preserving injection $\theta:B\to J$. Hence, $\im\psi\subseteq\phi(J)$, and for every $u\in B$, we have $a_{J,\phi}(\psi(u))=a_{J,\phi}(\phi(\theta(u)))=\kappa(\theta(u))=\kappa(u)$. Thus, $\ell(\psi)\subseteq x_{J,\phi}$.

Conversely, suppose that $\ell(\psi)\subseteq x_{J,\phi}$. Then $\theta=\phi^{-1}\circ\psi:B\to J$ is a well-defined injection, and the agreement of symbols gives $\kappa(\theta(u))=a_{J,\phi}(\psi(u))=\kappa(u)$ for $u\in B$. Therefore, $\theta$ is a template face with image in $J$, and $\phi\circ\theta=\psi$, so $v_{J,\phi}(\psi)=1$.
\end{proof}

Consequently, the molecule $v_{J,\phi}$ depends only on the active word and is constant on each orbit. Under \eqref{eq:orbit-identification}, it is the incidence vector of the labelled $t$-subsets of $x_{J,\phi}$, pulled back from orbits to individual role maps.

More generally, let $\cA$ be any family of separately indexed copies of $\Sigma^{\leq}$, called a $\Sigma^{\leq}$-family, and let $\gamma=(\gamma_{A,\theta})$ assign an integer to each $t$-face $\theta$ of each $A\in\cA$. For a $\Phi$-embedding $e=(A,\phi)$ of $A$, the molecule $v_e$ has entry $\gamma_{A,\theta}$ at $\phi\circ\theta$ and entry zero at faces without such a preimage. The \emph{decomposition lattice} is $\Lambda_\gamma(\Phi)=\langle v_e:e\in\cA(\Phi)\rangle_{\Z}$, where $\cA(\Phi)$ is the set of all such embeddings. An \emph{exact decomposition} of $U\in\Z^{\Phi_t}$ is a set $\mathcal D\subseteq\cA(\Phi)$ with $\sum_{e\in\mathcal D}v_e=U$. For the active templates in \eqref{eq:gamma}, we write $\Lambda(\Phi)=\Lambda_\gamma(\Phi)$, and we let $G=\one_{\Phi_t}\in\Z^{\Phi_t}$ be the all-one vector. By \eqref{eq:molecule-word}, the coordinate of $G$ at a role $t$-face represents the unit covering requirement for its labelled $t$-set. In \cite[Definition~2.21]{Keevash2018}, the molecule $v_e$ is written $\gamma(\phi)$, and $\Lambda_\gamma(\Phi)$ is written $\langle\gamma(\Phi)\rangle$.

Finally, we recall the remaining notions needed to state the decomposition theorem. For an orbit $O\in\Phi_t/\Sigma$ and $v\in\Z^{\Phi_t}$, let $v^O$ denote the vector that agrees with $v$ on $O$ and vanishes elsewhere.

\begin{definition}[see {\cite[Definitions~2.10, 2.22, 2.23, 2.25, 2.26, 2.36 and~2.37]{Keevash2018}}]\label{def:keevash}
Let $\cA$ be a $\Sigma^{\leq}$-family with an integer system $\gamma$ on the $t$-faces of its templates, let $\Phi$ be a $\Sigma$-adapted $R$-complex with $n=|V(\Phi)|$, and let $U\in\Z^{\Phi_t}$.
\begin{itemize}
\item The \emph{atoms} at $O\in\Phi_t/\Sigma$ are the distinct nonzero vectors $v_e^O$ with $e\in\cA(\Phi)$, and $\gamma[O]$ denotes their set. Repetitions of the same vector from different embeddings do not create different atoms. The system is \emph{elementary} if the atoms are linearly independent. We write $U\leq_\gamma V$ if $V-U$ is a sum of atoms with nonnegative integer coefficients, the empty sum included.
\item For $B\in\binom Rt$, let $\Sigma^B=\{\sigma|_{\sigma^{-1}(B)}:\sigma\in\Sigma\}$. The \emph{signature} of $(A,\theta)\in\cA_B$ is $\gamma^{A,\theta}=(\gamma_{A,\theta\circ\sigma})_{\sigma\in\Sigma^B}$, and its \emph{type} is its class under equality of signatures. Let $T_B$ be the set of types on $B$, write $\gamma^\tau$ for the common signature of $\tau\in T_B$, and call the class of the zero signature the zero type. For $\psi\in\Phi_B$ and $\tau\in T_B$, the \emph{orbit pattern} $\gamma[\psi]^\tau$ is the vector with $(\gamma[\psi]^\tau)_{\psi\circ\sigma}=(\gamma^\tau)_\sigma$ for $\sigma\in\Sigma^B$ and zero entries outside $\psi\Sigma$. For $e\in\cA(\Phi)$, we write $\tau_e(\psi)=\tau$ if $v_e^{\psi\Sigma}=\gamma[\psi]^\tau$.
\item Suppose that $\gamma$ is elementary and that $U$ lies in the integer span of the atoms. Each atom is supported on a single orbit, so $U^{\psi\Sigma}$ has a unique expansion in the atoms on $\psi\Sigma$. For a nonzero type $\tau\in T_B$ whose pattern $\gamma[\psi]^\tau$ is an atom, $U^\tau_\psi$ denotes the coefficient of that atom, and $U^\tau_\psi=0$ for the remaining nonzero types. Thus, $U^{\psi\Sigma}=\sum_{\tau\in T_B,\,\tau\neq0}U^\tau_\psi\,\gamma[\psi]^\tau$.
\item For $A\in\cA$, $B\in\binom Rt$ and $\psi\in\Phi_B$, the \emph{partial molecule} $v^{\rm part}_{A,\psi}$ has entries $(v^{\rm part}_{A,\psi})_{\psi\circ\sigma}=\gamma_{A,\sigma}$ for $\sigma\in\Sigma^B$, viewed as a $t$-face of $A$, and zero entries outside $\psi\Sigma$. Put $\gamma[U]_A=\{\psi\in\Phi_t:v^{\rm part}_{A,\psi}\leq_\gamma U\}$ and $\cA(\Phi,U)=\{e\in\cA(\Phi):v_e\leq_\gamma U\}$.
\item Let $0<\omega\leq1$ and $0<c<1$, and suppose that $\gamma$ is elementary and $U$ lies in the integer span of the atoms. Then $U$ is $(\gamma,c,\omega)$-\emph{regular} if there are real weights $(y_e)_{e\in\cA(\Phi,U)}$ with $\omega n^{t-L}\leq y_e\leq\omega^{-1}n^{t-L}$ such that, for every $B\in\binom Rt$, every $\psi\in\Phi_B$ and every nonzero type $\tau\in T_B$,
\[
(1-c)U_\psi^\tau\leq\sum_{e\in\cA(\Phi,U),\ \tau_e(\psi)=\tau}y_e
\leq(1+c)U_\psi^\tau.
\]
\item Let $\mathcal P=(\mathcal P_a:a\in I)$ be a family of subsets of $\Phi_t$. The pair $(\Phi,\mathcal P)$ is $(\omega,h)$-extendable if, for every rooted extension $(H,F,f)$ with $H\subseteq R(h)$ and every choice of pairwise disjoint sets $\mathcal H_a$ of $t$-faces of $H$ not in $H[F]$, at least $\omega|V(\Phi)|^{|V(H)\setminus F|}$ completions $g$ of $(H,F,f)$ satisfy $g\circ\theta\in\mathcal P_a$ for all $a\in I$ and $\theta\in\mathcal H_a$.
\end{itemize}
\end{definition}

The following result is the scalar case of Keevash's theorem, in which the parameters $q,r,D$ of \cite[Theorem~3.1]{Keevash2018} become $L,t,1$ and $R$ is identified with $[L]$.

\begin{lemma}[Keevash's theorem, see {\cite[Theorem~3.1]{Keevash2018}}]\label{lem:keevash}
Let $h=2^{50L^3}$ and $\delta=2^{-10^3L^5}$. There exist $n_0$ and $\omega_0>0$, depending only on $L$ and $t$, such that the following holds. Let $\cA$ be a $\Sigma^{\leq}$-family with an integer system $\gamma$ on the $t$-faces of its templates, let $\Phi$ be a $\Sigma$-adapted $R$-complex, and let $U\in\Z^{\Phi_t}$. Suppose that
\[
n=|V(\Phi)|>n_0,
\ \ n^{-\delta}<\omega<\omega_0,
\ \ 0<c\leq\omega^{h^{20}}.
\]
Assume that \textup{(1)} $\gamma$ is elementary, \textup{(2)} $U\in\Lambda_\gamma(\Phi)$, \textup{(3)} $U$ is $(\gamma,c,\omega)$-regular, and \textup{(4)} the pair $(\Phi,\gamma[U]_A)$ is $(\omega,h)$-extendable for every $A\in\cA$. Then $U$ has an exact decomposition.
\end{lemma}

To prove Theorem~\ref{thm:labelled}, we use the following consequence of Lemma~\ref{lem:keevash}, in which the target is the all-one vector $G$ and the fractional hypothesis is expressed directly as bounds on the total weight at each labelled $t$-set.

\begin{lemma}\label{lem:keevash-special}
Let $h=2^{50L^3}$ and $\delta=2^{-10^3L^5}$. There exist $\omega_0>0$ and $n_0$, depending only on $L$ and $t$, such that the following holds. Suppose that $n>n_0$, $n^{-\delta}<\omega<\omega_0$, and $0<c\leq\omega^{h^{20}}$. Let $\Phi$ be a $\Sigma$-adapted $R$-complex on $[n]$ satisfying \eqref{eq:full-skeleton}, and assume that
\begin{enumerate}
\item[(i)] $\Phi$ is $(\omega,h)$-extendable,
\item[(ii)] $G = \one_{\Phi_t} \in\Lambda(\Phi)$, and
\item[(iii)] there are real numbers $y_{J,\phi}$, one for every template embedding $(J,\phi)$, such that
\begin{equation}\label{eq:weight-range}
\omega n^{t-L}\leq y_{J,\phi}\leq\omega^{-1}n^{t-L}
\end{equation}
and, for every $\psi\in\Phi_t$,
\begin{equation}\label{eq:regularity}
1-c\leq\sum_{J,\phi}y_{J,\phi}v_{J,\phi}(\psi)\leq1+c.
\end{equation}
\end{enumerate}
Then there is a set $\mathcal D$ of template embeddings such that
\begin{equation}\label{eq:exact-decomposition}
\sum_{(J,\phi)\in\mathcal D}v_{J,\phi}=G.
\end{equation}
\end{lemma}

\begin{proof}
We apply Lemma~\ref{lem:keevash} with $U=G$ and the active system \eqref{eq:gamma}, choosing $\omega_0\leq1$, so that $c\leq\omega^{h^{20}}<1$. By Lemma~\ref{lem:embedding-identification}, the embeddings of $A_J$ are indexed by $\phi\in\Phi_R$, and their molecules are the vectors $v_{J,\phi}$. Thus, $\Lambda_\gamma(\Phi)=\Lambda(\Phi)$, and assumption~(ii) is condition~(2) of Lemma~\ref{lem:keevash}. It suffices to verify conditions (1), (3) and (4).

For (1), note that by \eqref{eq:molecule-word} each molecule is a $0$--$1$ vector that is constant on each orbit $O\in\Phi_t/\Sigma$. For any $\psi\in O$, the positive lower bound in \eqref{eq:regularity}, together with the positivity of the weights, implies that $v_e(\psi)=1$ for some embedding $e$. Thus,
\[
\gamma[O]=\{\one_O\}
\ \text{ and }\
G=\sum_{O\in\Phi_t/\Sigma}\one_O.
\]
The atoms are linearly independent, since evaluating a linear relation at a face in $O$ isolates the coefficient of $\one_O$. Hence, the system is elementary. Moreover, for each embedding $e$, let $\mathcal O_e=\{O\in\Phi_t/\Sigma:v_e^O=\one_O\}$. Then $G-v_e=\sum_{O\notin\mathcal O_e}\one_O$, which gives $v_e\leq_\gamma G$. Therefore, $\cA(\Phi,G)=\cA(\Phi)$.

For (3), fix $B\in\binom Rt$. Since every $\sigma\in\Sigma^B$ has image $B$, it follows from \eqref{eq:gamma} that $\gamma_{J,\theta\circ\sigma}=\one_{\im\theta\subseteq J}$ for all $\theta\in(A_J)_B$ and $\sigma\in\Sigma^B$. Hence, every signature is constant. The all-one signature occurs for $\theta=\id_B$ and any $J$ containing $B$, so there is exactly one nonzero type $\tau_+$ on $B$. For $\psi\in\Phi_B$, its orbit pattern and the coefficient of $G$ are
\[
\gamma[\psi]^{\tau_+}=\one_{\psi\Sigma}
\ \text{ and }\
G_\psi^{\tau_+}=1.
\]
Also, $\tau_e(\psi)=\tau_+$ exactly when $v_e(\psi)=1$. Since $\cA(\Phi,G)=\cA(\Phi)$ by the proof of (1), the weight sum for this type is
\[
\sum_{e\in\cA(\Phi,G),\ \tau_e(\psi)=\tau_+}y_e
=\sum_{e\in\cA(\Phi)}y_ev_e(\psi),
\]
where $y_e=y_{J,\phi}$ for $e=(J,\phi)$. Now \eqref{eq:weight-range} and \eqref{eq:regularity} give the $(\gamma,c,\omega)$-regularity of $G$.

Finally, for (4), let $B\in\binom Rt$ and $\psi\in\Phi_B$. The partial molecule for $A_J$ is $v^{\rm part}_{J,\psi}=\one_{B\subseteq J}\one_{\psi\Sigma}$, since its template entries are $\gamma_{J,\sigma}$ with $\im\sigma=B$. Subtracting this vector from $G$ leaves a sum of atoms with nonnegative coefficients, so $\gamma[G]_{A_J}=\Phi_t$ for every $J$. Consequently, the restrictions in the definition of extendability for the pair $(\Phi,\gamma[G]_{A_J})$ only require the selected $t$-faces to be mapped into $\Phi_t$, which holds for every $\Phi$-embedding. Hence, assumption~(i) gives condition~(4).

All conditions of Lemma~\ref{lem:keevash} hold for $U=G$. The resulting exact decomposition is a set $\mathcal D\subseteq\cA(\Phi)$ satisfying \eqref{eq:exact-decomposition}, which completes the proof.
\end{proof}

\section{Proof of Theorem \texorpdfstring{\ref*{thm:labelled}}{1.5} via Keevash's theorem}\label{sec:proof}

In this section, we prove Theorem~\ref{thm:labelled} via Lemma~\ref{lem:keevash-special} and then deduce Theorem~\ref{thm:main}.

\paragraph{Proof overview.} As outlined in Section~\ref{subsec:approach}, we prescribe one random labelling on each coordinate $(t+1)$-set and permit only those role maps whose labels agree with these prescriptions. Two permitted blocks sharing a coordinate $(t+1)$-set then carry the same labels on it, and hence share a labelled $t$-set. Consequently, an exact labelled cover by permitted blocks automatically has support intersections of size at most $t$, see Lemma~\ref{lem:extraction}. The existence of such an exact cover follows from Lemma~\ref{lem:keevash-special}, whose three hypotheses are extendability, integer feasibility and fractional regularity. For the random host, Lemma~\ref{lem:extendable-random} establishes extendability together with uniform counts of full extensions, Lemma~\ref{lem:host-lattice} gives integer feasibility, and Lemma~\ref{lem:regularity} derives fractional regularity from these counts. The proofs of these three lemmas are given in Section~\ref{sec:lemmas}.

We now define the random host. For each coordinate set $S\in\binom{[n]}{t+1}$, choose a map $\chi_S:S\to[m]$ uniformly from the $m^{t+1}$ possible maps, independently for different sets $S$. We call $\chi_S$ the \emph{selector} of $S$ and write $\boldsymbol\chi=(\chi_S)_{S\in\binom{[n]}{t+1}}$. For $|B|\leq t$, the set $B$ has no $(t+1)$-subset, and we set $\Phi_B=\Phi_B^0=\Inj(B,[n])$. For $|B|\geq t+1$, let $\Phi_B$ consist of all injections $\phi:B\to[n]$ satisfying
\begin{equation}\label{eq:selector}
\chi_{\phi(C)}(\phi(u))=\kappa(u)
\ \text{ for all } C\in\binom{B}{t+1} \text{ and } u\in C.
\end{equation}
For each $C$, the image $S=\phi(C)$ determines which selector is used, so any two maps that send $(t+1)$-subsets of their domains onto the same set $S$ are checked against the same selector $\chi_S$. We write $\Phi=\Phi(\boldsymbol\chi)=(\Phi_B:B\subseteq R)$ and call it the \emph{random selector complex}. By construction, $\Phi$ and $\Phi^0$ have the same faces of size at most $t$.

\begin{lemma}\label{lem:random-adapted}
For every realization of the selectors $\boldsymbol\chi$, the complex $\Phi=\Phi(\boldsymbol\chi)$ is a $\Sigma$-adapted subcomplex of $\Phi^0$ satisfying \eqref{eq:full-skeleton}, and $V(\Phi)=[n]$.
\end{lemma}

\begin{proof}
For every $B\subseteq R$, we have $\Phi_B\subseteq\Inj(B,[n])$. If $\phi\in\Phi_B$ and $B'\subseteq B$, then every $(t+1)$-subset of $B'$ is also a $(t+1)$-subset of $B$, so all selector conditions required for $\phi|_{B'}$ are already satisfied by $\phi$. Hence, $\phi|_{B'}\in\Phi_{B'}$, and $\Phi$ is a subcomplex of $\Phi^0$. By construction, $\Phi$ satisfies \eqref{eq:full-skeleton}, and since $t\geq1$, every coordinate occurs in a singleton face, so $V(\Phi)=[n]$.

To prove adaptedness, it suffices to consider $B$ with $|B|\geq t+1$. Let $\sigma\in\Sigma$ and $\phi\in\Phi_B$, and recall that $\phi\sigma$ has domain $\sigma^{-1}(B)$. For any $C\in\binom{\sigma^{-1}(B)}{t+1}$ and $u\in C$, we have $\sigma(C)\in\binom{B}{t+1}$ and
\[
\chi_{(\phi\sigma)(C)}\bigl((\phi\sigma)(u)\bigr)
=\chi_{\phi(\sigma(C))}\bigl(\phi(\sigma(u))\bigr)
=\kappa(\sigma(u))=\kappa(u),
\]
where the second equality uses $\phi\in\Phi_B$ and the last one uses the color preservation of $\sigma$. Thus, $\phi\sigma\in\Phi_{\sigma^{-1}(B)}$.
\end{proof}

We next define the quantities appearing in the extension counts and the fractional weights. For $n\geq L$, put
\begin{equation}\label{eq:p-rho}
p=m^{-(t+1)},
\ \ \rho=p^{\binom L{t+1}}
\ \text{ and }\
D_n=\rho\cdot\fall{n-t}{L-t},
\end{equation}
where $\fall xk=\prod_{i=0}^{k-1}(x-i)$ is the falling factorial, with $\fall x0=1$. For each $J\in\binom Rw$, write $j_a=|J\cap P_a|$ for $a\in[m]$, put $M(J)=\prod_{a=1}^m\fall w{j_a}$, and define
\begin{equation}\label{eq:alpha}
\alpha_J=\frac{(w-t)!}{m^{w-t}M(J)}.
\end{equation}
Since $0\leq j_a\leq w$, we have $\alpha_J>0$, and $\alpha_J$ depends only on the color counts $(j_a)_{a\in[m]}$. To each template embedding $(J,\phi)$, we assign the weight
\begin{equation}\label{eq:weights}
y_{J,\phi}=\frac{\alpha_J}{D_n},
\end{equation}
and we set
\[
a_0=\frac{1}{\rho}\min_{J\in\binom Rw}\alpha_J
\ \text{ and }\
a_1=\frac{2}{\rho}\max_{J\in\binom Rw}\alpha_J.
\]
The constants $\rho$, $a_0$ and $a_1$ depend only on $m,w,t$.

The next three lemmas show that the random selector complex satisfies the hypotheses of Lemma~\ref{lem:keevash-special} with positive probability. Their proofs are given in Section~\ref{sec:lemmas}.

\begin{lemma}\label{lem:extendable-random}
For every fixed $h\in\N$, there is a constant $\omega_1=\omega_1(m,w,t,h)>0$ such that, for every fixed $c>0$, with probability $1-o(1)$ as $n\to\infty$, the random selector complex $\Phi$ is $(\omega_1,h)$-extendable and satisfies, with $D_n$ as in \eqref{eq:p-rho},
\begin{equation}\label{eq:full-extension-count}
(1-c)D_n
\leq\bigl|\{\phi\in\Phi_R:\phi|_B=\psi\}\bigr|
\leq(1+c)D_n
\ \text{ for all } B\in\binom Rt \text{ and } \psi\in\Phi_B.
\end{equation}
\end{lemma}

\begin{lemma}\label{lem:host-lattice}
Fix $\eta>0$. For all sufficiently large $n$ satisfying \eqref{eq:admissible}, if $\Phi$ is a $\Sigma$-adapted $(\eta,3t^2)$-extendable $R$-complex on $[n]$ satisfying \eqref{eq:full-skeleton}, then $G= \one_{\Phi_t}\in\Lambda(\Phi)$.
\end{lemma}

\begin{lemma}\label{lem:regularity}
If an $R$-complex $\Phi$ on $[n]$ satisfies \eqref{eq:full-skeleton} and \eqref{eq:full-extension-count} with error $c$, then the weights in \eqref{eq:weights} satisfy \eqref{eq:regularity}. Moreover, for all sufficiently large $n$ and every template embedding $(J,\phi)$, we have $a_0n^{t-L}\leq y_{J,\phi}\leq a_1n^{t-L}$.
\end{lemma}

The next lemma shows that every exact decomposition in a selector complex defines a family as in Theorem~\ref{thm:labelled}.

\begin{lemma}\label{lem:extraction}
Let $\Phi=\Phi(\boldsymbol\chi)$ be the selector complex of some realization $\boldsymbol\chi$, and let $\mathcal D\subseteq\cA(\Phi)$ satisfy \eqref{eq:exact-decomposition}. For $e=(J,\phi)\in\mathcal D$, let $x_e=x_{J,\phi}$ be its active word \eqref{eq:active-word}. Then the family $\mathcal B=\{x_e:e\in\mathcal D\}\subseteq\cT_w(n,m)$ satisfies properties (i) and (ii) of Theorem~\ref{thm:labelled}.
\end{lemma}

\begin{proof}
Each $x_e$ is a labelled $w$-set, since $\phi$ is injective and $|J|=w$. To prove (i), fix a labelled $t$-set $u\in\cT_t(n,m)$. By Lemmas~\ref{lem:random-adapted} and~\ref{lem:orbit-identification}, there is a role $t$-face $\psi\in\Phi_t$ with $\ell(\psi)=u$. It then follows from \eqref{eq:exact-decomposition} and \eqref{eq:molecule-word} that
\[
1=G(\psi)=\sum_{e\in\mathcal D}v_e(\psi)
=\sum_{e\in\mathcal D}\one_{u\subseteq x_e}.
\]
Every summand is $0$ or $1$, so exactly one selected embedding $e$ has an active word containing $u$. If two distinct selected embeddings $e$ and $e'$ had $x_e=x_{e'}=W$, then every labelled $t$-subset of $W$ would be covered by both, contrary to this uniqueness. Thus, the map $e\mapsto x_e$ is injective on $\mathcal D$, and every labelled $t$-set lies in exactly one member of $\mathcal B$.

To prove (ii), suppose for a contradiction that two distinct members $x=x_{J,\phi}$ and $y=x_{J',\phi'}$ of $\mathcal B$ have supports meeting in at least $t+1$ coordinates, and choose $S\subseteq\supp(x)\cap\supp(y)$ with $|S|=t+1$. Since $\supp(x)=\phi(J)$, the preimage $C=\phi^{-1}(S)$ is a $(t+1)$-subset of $J$. The full map $\phi\in\Phi_R$ satisfies the selector condition \eqref{eq:selector} on $C$, so $\chi_S(\phi(u))=\kappa(u)$ for $u\in C$. By \eqref{eq:active-word}, $\kappa(u)$ is the symbol of $x$ at the coordinate $\phi(u)$, and hence $x|_S=\chi_S$. Applying the same argument to $\phi'$, we obtain $y|_S=\chi_S$. For any $T\in\binom St$, both words then contain the labelled $t$-set $(T,\chi_S|_T)$, which contradicts (i).
\end{proof}

\begin{proof}[Proof of Theorem~\ref{thm:labelled}]
Let $n_0,h,\delta,\omega_0$ be the constants in Lemma~\ref{lem:keevash-special}, let $\omega_1$ be the constant given by Lemma~\ref{lem:extendable-random} for this $h$, and let $a_0,a_1$ be the constants defined after \eqref{eq:weights}. Choose $\omega>0$ with
\[
\omega<\min\{\omega_0,1\},
\ \ \omega\leq\omega_1,
\ \ \omega\leq a_0
\ \text{ and }\
\omega\leq a_1^{-1},
\]
and then set $c=\omega^{h^{20}}$. All these constants depend only on $m,w,t$. By Lemma~\ref{lem:extendable-random}, the constant $\omega_1$ does not depend on $c$. We can therefore choose $\omega$ first and then impose the much smaller error $c=\omega^{h^{20}}$ without changing $\omega_1$.

Let $n$ satisfy \eqref{eq:admissible} and be sufficiently large. In particular, we require $n>n_0$, $n^{-\delta}<\omega$ and $n\geq L$, together with the thresholds in Lemmas~\ref{lem:extendable-random}, \ref{lem:host-lattice} and~\ref{lem:regularity} for the fixed parameters just chosen. These requirements impose only a threshold depending on $m,w,t$. By Lemma~\ref{lem:extendable-random}, with probability $1-o(1)$, the random selector complex is $(\omega_1,h)$-extendable and satisfies \eqref{eq:full-extension-count} with error $c$. For sufficiently large $n$, this probability is positive, and we fix a realization $\boldsymbol\chi$ for which both properties hold. Let $\Phi=\Phi(\boldsymbol\chi)$. By Lemma~\ref{lem:random-adapted}, $\Phi$ is a $\Sigma$-adapted $R$-complex on $[n]$ satisfying \eqref{eq:full-skeleton}. We now verify the three hypotheses of Lemma~\ref{lem:keevash-special}.

Since $\omega\leq\omega_1$, it follows from Definition~\ref{def:extendability} that $\Phi$ is $(\omega,h)$-extendable, which is hypothesis~(i). As $h\geq3t^2$, the complex $\Phi$ is also $(\omega,3t^2)$-extendable, and applying Lemma~\ref{lem:host-lattice} with $\eta=\omega$ yields $G\in\Lambda(\Phi)$, which is hypothesis~(ii). For hypothesis~(iii), we use the weights $y_{J,\phi}$ in \eqref{eq:weights}. By Lemma~\ref{lem:regularity}, they satisfy \eqref{eq:regularity} and $a_0n^{t-L}\leq y_{J,\phi}\leq a_1n^{t-L}$, so the choice of $\omega$ gives \eqref{eq:weight-range}.

Lemma~\ref{lem:keevash-special} now provides a set $\mathcal D\subseteq\cA(\Phi)$ satisfying \eqref{eq:exact-decomposition}, and by Lemma~\ref{lem:extraction}, the family $\mathcal B=\{x_e:e\in\mathcal D\}$ has properties (i) and (ii). Since all thresholds and constants above depend only on $m,w,t$, this holds for every sufficiently large $n$ satisfying \eqref{eq:admissible}.
\end{proof}

\begin{proof}[Proof of Theorem~\ref{thm:main}]
Put $m=q-1\geq2$. For every sufficiently large $n$ satisfying \eqref{eq:admissible}, Theorem~\ref{thm:labelled} supplies a family $\mathcal B\subseteq\cT_w(n,m)$ with exact labelled coverage and support intersections of size at most $t$. By Lemma~\ref{lem:distance-transfer}, the associated code is a $GS(t,w,n,q)$.
\end{proof}

\section{Proofs of Lemmas \texorpdfstring{\ref*{lem:extendable-random}, \ref*{lem:host-lattice} and \ref*{lem:regularity}}{3.2, 3.3 and 3.4}}\label{sec:lemmas}

\subsection{Proof of Lemma \texorpdfstring{\ref*{lem:extendable-random}}{3.2}}\label{subsec:extendability}

In this subsection, we prove Lemma~\ref{lem:extendable-random}. We first count the possible images of one additional vertex, directly in terms of coordinates and selectors. Suppose that $U\subseteq[n]$ is a bounded set of coordinates already in use, and that $b:U\to[m]$ records their assigned symbols in a partial embedding. A new vertex with prescribed symbol $a\in[m]$ must be mapped to some $x\in[n]\setminus U$, and a family $\mathcal S\subseteq\binom Ut$ specifies which selectors on the sets $T\cup\{x\}$ must be checked. We require each such selector to agree with $b$ on $T$ and to assign the symbol $a$ to $x$.

\begin{lemma}\label{lem:one-vertex}
Fix an integer $N\geq1$ and a real number $0<\varepsilon<1$. With probability tending to one as $n\to\infty$, the following holds simultaneously for every set $U\subseteq[n]$ with $|U|\leq N$, every map $b:U\to[m]$, every $a\in[m]$ and every family $\mathcal S\subseteq\binom Ut$. The number of coordinates $x\in[n]\setminus U$ satisfying
\begin{equation}\label{eq:one-vertex-condition}
\chi_{T\cup\{x\}}|_T=b|_T
\ \text{ and }\
\chi_{T\cup\{x\}}(x)=a
\ \text{ for all } T\in\mathcal S
\end{equation}
lies in the interval
\begin{equation}\label{eq:one-vertex-count}
\bigl[(1-\varepsilon)p^{|\mathcal S|}(n-|U|),\,(1+\varepsilon)p^{|\mathcal S|}(n-|U|)\bigr].
\end{equation}
\end{lemma}

\begin{proof}
Fix a tuple $Q=(U,b,a,\mathcal S)$ as in the statement. All randomness comes from the selectors, which are independent and uniformly distributed. For $x\in[n]\setminus U$ and $T\in\mathcal S$, let $E_T(x)$ be the event that $\chi_{T\cup\{x\}}|_T=b|_T$ and $\chi_{T\cup\{x\}}(x)=a$, and let $I_x$ be the indicator of the event that $E_T(x)$ holds for every $T\in\mathcal S$. Thus, the number of admissible coordinates for $Q$ is $Z_Q=\sum_{x\in[n]\setminus U}I_x$.

Since $|T|=t$ and $x\notin U$, the set $T\cup\{x\}$ has size $t+1$, and $E_T(x)$ prescribes its selector at every point. Exactly one of the $m^{t+1}$ maps $T\cup\{x\}\to[m]$ satisfies these requirements, so $\Prob(E_T(x))=m^{-(t+1)}=p$. If $T,T'\in\mathcal S$ are distinct, then $T\cup\{x\}\neq T'\cup\{x\}$, because removing $x$ recovers $T$ and $T'$. Hence, the events $E_T(x)$ with $T\in\mathcal S$ depend on distinct selectors and are mutually independent, even when the sets $T$ intersect. Consequently,
\[
\Prob(I_x=1)=\Prob\left(\bigcap_{T\in\mathcal S}E_T(x)\right)
=\prod_{T\in\mathcal S}\Prob(E_T(x))=p^{|\mathcal S|},
\]
where the empty product is $1$. Moreover, every set $T\cup\{x\}$ with $T\subseteq U$ has $x$ as its only point outside $U$. Thus, the indicators $I_x$ with $x\in[n]\setminus U$ depend on pairwise disjoint collections of the independent selectors, and hence they are mutually independent. It follows that $Z_Q$ has the binomial distribution $\operatorname{Bin}(n-|U|,p^{|\mathcal S|})$, with expectation $\mu_Q=\E Z_Q=p^{|\mathcal S|}(n-|U|)$.

Let $F_Q$ be the event that $|Z_Q-\mu_Q|>\varepsilon\mu_Q$. If $\mathcal S=\varnothing$, then $Z_Q=n-|U|=\mu_Q$ deterministically, and $\Prob(F_Q)=0$. In particular, if $N<t$, then all families $\mathcal S$ are empty, and the lemma holds deterministically. In general, we have $|\mathcal S|\leq\binom{|U|}t\leq\binom Nt$, and since $0<p<1$, for $n\geq2N$ we obtain
\[
\mu_Q=p^{|\mathcal S|}(n-|U|)\geq p^{\binom Nt}(n-N)
\geq\frac12np^{\binom Nt}.
\]
The multiplicative Chernoff bound, applied with $0<\varepsilon<1$, gives
\[
\Prob(F_Q)\leq2\exp\Bigl(-\frac{\varepsilon^2\mu_Q}{3}\Bigr)\leq2e^{-c_1n},
\ \text{ where }\
c_1=\frac{\varepsilon^2}{6}p^{\binom Nt}>0.
\]
This bound also holds when $\mathcal S=\varnothing$, so it holds for every tuple $Q$ with the same constant $c_1$.

Let $\mathcal Q_n$ be the set of all tuples $Q=(U,b,a,\mathcal S)$ as in the statement. For $n\geq N$,
\[
|\mathcal Q_n|=\sum_{j=0}^{N}\binom njm^jm\,2^{\binom jt}
\leq m\,2^{\binom Nt}\sum_{j=0}^{N}\binom njm^j\leq Cn^N,
\]
where $C>0$ depends only on $m,t,N$. By the union bound, the probability that $F_Q$ occurs for at least one $Q\in\mathcal Q_n$ is at most $|\mathcal Q_n|\cdot2e^{-c_1n}\leq2Cn^Ne^{-c_1n}=o(1)$, since $N$ and $c_1>0$ are fixed. Taking complements, with probability $1-o(1)$ we have $|Z_Q-\mu_Q|\leq\varepsilon\mu_Q$ for every $Q\in\mathcal Q_n$, which is \eqref{eq:one-vertex-count} for all tuples simultaneously. Once a realization of the selectors in this event is fixed, the estimate holds for every tuple, and hence it remains valid for tuples chosen adaptively in the embedding arguments below.
\end{proof}

\begin{proof}[Proof of Lemma~\ref{lem:extendable-random}]
We first prove the extendability. Fix $h$ and put $N=Lh$. By Lemma~\ref{lem:one-vertex} with $\varepsilon=1/2$, with probability $1-o(1)$, its conclusion holds for every tuple $(U,b,a,\mathcal S)$ with $|U|\leq N$. Fix a realization $\boldsymbol\chi$ in this event, and let $\Phi=\Phi(\boldsymbol\chi)$. Let $H\subseteq R(h)$, let $F\subseteq V(H)$, and let $f:F\to[n]$ be a $\Phi$-embedding of $H[F]$. Since $H\subseteq R(h)$, we have $|V(H)|\leq Lh=N$. Each abstract vertex $y=(u,j)\in V(H)$ has the fixed color $\kappa_h(y)=\kappa(u)$, which prescribes the symbol at its image $g(y)$ under an embedding $g$. Let
\[
\cY(H)=\left\{Y\in\binom{V(H)}{t+1}:Y\subseteq\im\theta
\text{ for some face }\theta\text{ of }H\right\}.
\]
A face $\theta\in H_B$ sends each role $u\in B$ to an abstract vertex of color $\kappa(u)$, and $C\mapsto\theta(C)$ is a bijection from $\binom B{t+1}$ to the $(t+1)$-subsets of $\im\theta$. By \eqref{eq:selector}, an injective map $g:V(H)\to[n]$ is therefore a $\Phi$-embedding of $H$ if and only if
\begin{equation}\label{eq:Y-condition}
\chi_{g(Y)}(g(y))=\kappa_h(y)
\ \text{ for every } Y\in\cY(H) \text{ and every } y\in Y.
\end{equation}
Since $g$ is injective, $g(Y)=g(Y')$ implies $Y=Y'$, so distinct sets in $\cY(H)$ query distinct selectors, and a set occurring in several faces receives the same prescription $i\mapsto\kappa_h(g^{-1}(i))$ each time. Thus, these requirements are consistent. Moreover, the conditions in \eqref{eq:Y-condition} with $Y\subseteq F$ are already satisfied by $f$, since restricting a face whose image contains $Y$ to the preimage of $Y$ gives a face of $H[F]$, and $f$ maps it into $\Phi$.

It remains to embed the vertices outside $F$. Order them as $y_1,\ldots,y_v$, where $v=|V(H)\setminus F|$, and put $F_k=F\cup\{y_1,\ldots,y_k\}$ for $0\leq k\leq v$. Suppose that an injective map $g_{k-1}:F_{k-1}\to[n]$ extending $f$ has been constructed and satisfies all conditions in \eqref{eq:Y-condition} whose sets lie in $F_{k-1}$. Let $U_k=g_{k-1}(F_{k-1})$, and define $b_k:U_k\to[m]$ by $b_k(i)=\kappa_h(g_{k-1}^{-1}(i))$, which is well defined because $g_{k-1}$ is injective. The new conditions are those indexed by the sets $Y\in\cY(H)$ with $y_k\in Y\subseteq F_k$. All vertices of $Y\setminus\{y_k\}$ are already embedded, and we put
\[
\mathcal S_k=\bigl\{g_{k-1}(Y\setminus\{y_k\}):Y\in\cY(H),\ y_k\in Y\subseteq F_k\bigr\}
\subseteq\binom{U_k}{t}.
\]
By the injectivity of $g_{k-1}$, distinct sets $Y$ give distinct members of $\mathcal S_k$. The requirements on a proposed image $x\in[n]\setminus U_k$ of $y_k$ are exactly the conditions \eqref{eq:one-vertex-condition} with $U=U_k$, $b=b_k$, $a=\kappa_h(y_k)$ and $\mathcal S=\mathcal S_k$. Since $|U_k|\leq N$ and $|\mathcal S_k|\leq\binom Nt$, it follows from Lemma~\ref{lem:one-vertex} that, for $n\geq2N$, the number of admissible choices is at least
\[
\frac12p^{|\mathcal S_k|}(n-|U_k|)\geq\frac12p^{\binom Nt}(n-N)
\geq\frac14p^{\binom Nt}n=:\beta n.
\]
Choosing such an $x$ extends $g_{k-1}$ to an injective map $g_k$ on $F_k$ satisfying every condition inside $F_k$.

Every set $Y\in\cY(H)$ either lies in $F$ or has a unique last vertex in the order $y_1,\ldots,y_v$, and its condition is checked when that vertex is embedded. Consequently, every sequence of admissible choices yields a completion of $(H,F,f)$, and distinct sequences yield distinct completions. Hence, the number of completions is at least $(\beta n)^v\geq\beta^Nn^v$, because $0<\beta\leq1$ and $v\leq N$. This also holds when $v=0$, in which case $f$ itself is the unique completion. Thus, $\Phi$ is $(\omega_1,h)$-extendable with
\[
\omega_1=\beta^N,
\ \ \beta=\frac14p^{\binom Nt}
\ \text{ and }\
N=mwh,
\]
which depend only on $m,w,t,h$. Since the argument applies to every rooted extension on the same event, no further union bound is needed.

It remains to prove \eqref{eq:full-extension-count}. Fix $c>0$ and put $c_*=\min\{c,1/2\}$. It suffices to prove the estimate with $c_*$ in place of $c$, since this gives a stronger bound. By continuity at $\varepsilon=0$, we can choose $0<\varepsilon<1$ such that $(1-\varepsilon)^{L-t}\geq1-c_*$ and $(1+\varepsilon)^{L-t}\leq1+c_*$. By Lemma~\ref{lem:one-vertex}, with probability $1-o(1)$, its conclusion holds for $N=L$ and this $\varepsilon$. Fix a realization of the selectors in this event.

Let $B\in\binom Rt$ and $\psi\in\Phi_B$. Order the remaining roles as $R\setminus B=\{u_{t+1},\ldots,u_L\}$, put $B_t=B$ and $B_j=B\cup\{u_{t+1},\ldots,u_j\}$ for $t<j\leq L$, and let
\[
\mathcal E_j=\{\phi_j\in\Phi_{B_j}:\phi_j|_B=\psi\}.
\]
Thus, $\mathcal E_t=\{\psi\}$, and the desired count is $|\mathcal E_L|$. Fix $\phi_j\in\mathcal E_j$ with $t\leq j<L$. To extend it to $B_{j+1}$, we choose an unused coordinate $x$ for the role $u_{j+1}$. All selector conditions on subsets of $B_j$ are already satisfied, and the remaining ones are those on the sets $C\cup\{u_{j+1}\}$ with $C\in\binom{B_j}t$. Put $U=\phi_j(B_j)$, define $b$ by $b(\phi_j(u))=\kappa(u)$ for $u\in B_j$, and let $a=\kappa(u_{j+1})$. Since $\phi_j$ is a bijection from $B_j$ to $U$, the family of coordinate sets to be tested is
\[
\mathcal S=\left\{\phi_j(C):C\in\binom{B_j}{t}\right\}=\binom Ut,
\]
so $|U|=j$ and $|\mathcal S|=\binom jt$. By Lemma~\ref{lem:one-vertex}, the number of extensions of $\phi_j$ lies between $(1-\varepsilon)p^{\binom jt}(n-j)$ and $(1+\varepsilon)p^{\binom jt}(n-j)$. These bounds are uniform over all $\phi_j\in\mathcal E_j$, and every member of $\mathcal E_{j+1}$ has a unique restriction in $\mathcal E_j$. Summing over $\phi_j\in\mathcal E_j$, we obtain
\[
|\mathcal E_j|(1-\varepsilon)p^{\binom jt}(n-j)
\leq|\mathcal E_{j+1}|
\leq|\mathcal E_j|(1+\varepsilon)p^{\binom jt}(n-j).
\]
Starting from $|\mathcal E_t|=1$ and using
\[
\prod_{j=t}^{L-1}(n-j)=\fall{n-t}{L-t}
\ \text{ and }\
\sum_{j=t}^{L-1}\binom jt=\binom L{t+1},
\]
we obtain $(1-\varepsilon)^{L-t} \cdot \rho \cdot \fall{n-t}{L-t}\leq|\mathcal E_L|\leq(1+\varepsilon)^{L-t} \cdot \rho \cdot \fall{n-t}{L-t}$. By the choice of $\varepsilon$, this gives $(1-c_*)D_n\leq|\mathcal E_L|\leq(1+c_*)D_n$, and since $c_*\leq c$, it proves \eqref{eq:full-extension-count}. The one-vertex estimate holds for all tuples simultaneously, so this count holds for every $B$ and $\psi$ at once. The stepwise bounds are multiplied through the displayed recurrence for $|\mathcal E_j|$ on the fixed realization, so no independence between the steps is required.

Finally, for fixed $h$ and $c$, each of the two events used above has failure probability $o(1)$, so both conclusions hold simultaneously with probability $1-o(1)$. The constant $\omega_1$ was obtained with the fixed tolerance $1/2$ and does not depend on $c$.
\end{proof}

\subsection{Proof of Lemma \texorpdfstring{\ref*{lem:host-lattice}}{3.3}}\label{subsec:lattice}

In this subsection, we prove Lemma~\ref{lem:host-lattice}. We first establish $G\in\Lambda(\Phi^0)$ by viewing labelled sets as subsets of the coordinate--symbol space $X=[n]\times[m]$, where the covering problem becomes an ordinary hypergraph lattice problem. We then use the common skeleton to transfer this membership to $\Lambda(\Phi)$.

Define a $[w]$-complex $\Xi$ on $X$ by
\[
\Xi_B=\{f\in\Inj(B,X):\im f\text{ is transversal}\}
\ \text{ for } B\subseteq[w].
\]
It is closed under restriction because every subset of a transversal set is transversal, and it is $\Sym([w])$-adapted in the sense of \eqref{eq:adapted}, since relabelling the domain of an injection does not change its image. Every singleton is transversal, so $V(\Xi)=X$ and $|V(\Xi)|=mn$. For $j\leq w$, the images of the faces of size $j$ are exactly the transversal $j$-subsets of $X$, which correspond to $\cT_j(n,m)$ under \eqref{eq:label-transversal}.

\begin{lemma}\label{lem:transversal-extendable}
The complex $\Xi$ is $(2^{-3wt^2},3t^2)$-extendable whenever $n\geq6wt^2$.
\end{lemma}

\begin{proof}
Let $H\subseteq[w](3t^2)$ be a subcomplex, let $F\subseteq V(H)$, and let $f:F\to X$ be a $\Xi$-embedding of $H[F]$. Since $H$ has at most $3wt^2$ vertices, we have $v=|V(H)\setminus F|\leq3wt^2$. We embed the vertices of $V(H)\setminus F$ one at a time, choosing for each of them an image $(i,a)\in X$ whose first coordinate $i$ is used neither by $f(F)$ nor by a previously embedded vertex. At most $3wt^2$ first coordinates are forbidden at each step, so there are at least $m(n-3wt^2)\geq mn/2=|X|/2$ choices for each image. This gives an injective map $g:V(H)\to X$ extending $f$.

We check that $g$ is a $\Xi$-embedding. Let $\theta\in H_B$, and put $B_F=\{u\in B:\theta(u)\in F\}$. By restriction closure, $\theta|_{B_F}$ is a face of $H[F]$, so $f(\theta(B_F))$ is transversal. The other vertices of $\im\theta$ were assigned first coordinates that are distinct from all root first coordinates and from one another. Hence, $g(\im\theta)$ is transversal, and $g\circ\theta\in\Xi_B$. Note that the argument does not require $f(F)$ to be transversal. Different root vertices may share a first coordinate if they do not occur together in a face, and only the root part of each individual face needs to be transversal.

Every sequence of choices produces a different completion. Hence, the number of completions is at least $(|X|/2)^v=2^{-v}|X|^v\geq2^{-3wt^2}|X|^v$, which also holds for $v=0$.
\end{proof}

We next recall the lattice criterion for hypergraph decompositions.

\begin{definition}[Hypergraph degrees and divisibility]\label{def:hypergraph-divisibility}
An $r$-graph on a set $V$ is a family of $r$-subsets of $V$. For an $r$-graph $K$ on $[k]$ and $0\leq i\leq r$, set
\[
\deg_K(I)=|\{e\in K:I\subseteq e\}|
\ \text{ for } I\in\binom{[k]}{i},
\ \text{ and }\
\gcd\nolimits_i(K)=\gcd\left\{\deg_K(I):I\in\binom{[k]}{i}\right\}.
\]
Let $U$ be a nonnegative integer vector indexed by a family of $r$-subsets of a set $V$, extended by zero to $\binom Vr$, and put $\deg_U(S)=\sum_{e\in\binom Vr,\,S\subseteq e}U(e)$ for $S\subseteq V$ with $|S|\leq r$. We say that $U$ is $K$-\emph{divisible} if $\gcd_i(K)\mid\deg_U(S)$ for every $0\leq i\leq r$ and every $S\in\binom Vi$, where a condition with $\gcd_i(K)=0$ means $\deg_U(S)=0$. The case $S=\varnothing$ is the condition on the total weight of $U$.
\end{definition}

\begin{lemma}[see {\cite[Lemma~5.20]{Keevash2018}}]\label{lem:hypergraph-lattice}
Let $K$ be an $r$-graph on $[k]$, and let $\Psi$ be a $\Sym([k])$-adapted, $(\eta,3r^2)$-extendable $[k]$-complex. Suppose that $N_*=|V(\Psi)|$ is sufficiently large and $\eta>N_*^{-1/2}$. Write $\im\Psi_r=\{\im\psi:\psi\in\Psi_r\}$, and for each full embedding $\phi\in\Psi_{[k]}$ let $\phi(K)=\{\phi(e):e\in K\}\subseteq\im\Psi_r$. Then a nonnegative integer vector $U$ on $\im\Psi_r$ is $K$-divisible if and only if there are integers $c_\phi$ such that $U=\sum_{\phi\in\Psi_{[k]}}c_\phi\,\one_{\phi(K)}$.
\end{lemma}

This is \cite[Lemma~5.20]{Keevash2018} applied with $(q,H,\Phi,\omega,n,G)=(k,K,\Psi,\eta,N_*,U)$ and $s=3r^2$. The ground set $\Phi^\circ_r$ in that lemma is $\im\Psi_r$, and the copy family $H(\Phi)$ consists of the edge sets $\phi(K)$ with $\phi\in\Psi_{[k]}$. Indexing the generators by embeddings leaves the integer span unchanged, since embeddings with the same edge set give the same indicator vector.

Let $K=\binom{[w]}{t}$ be the complete $t$-graph on $[w]$. The images $\im\Xi_t$ of the $t$-faces of $\Xi$ are exactly the transversal $t$-subsets of $X$, and we take the target
\[
G_X=\one_{\im\Xi_t}\in\Z^{\im\Xi_t},
\]
which assigns unit weight to $\widehat u$ for every $u\in\cT_t(n,m)$.

\begin{lemma}\label{lem:transversal-divisible}
If $n$ satisfies \eqref{eq:admissible}, then $G_X$ is $K$-divisible.
\end{lemma}

\begin{proof}
For $0\leq i\leq t$, each $i$-subset of $[w]$ lies in $\binom{w-i}{t-i}$ edges of $K$, so $\gcd_i(K)=\binom{w-i}{t-i}$. The degree $\deg_{G_X}(I)$ counts the transversal $t$-sets containing $I$. If $I$ is a transversal $i$-set, such a $t$-set is obtained by choosing $t-i$ of the $n-i$ unused coordinates and one of the $m$ symbols for each of them. If $I$ is not transversal, no set containing it is transversal. Thus,
\[
\deg_{G_X}(I)=
\begin{cases}
m^{t-i}\binom{n-i}{t-i},&\text{if $I$ is transversal},\\
0,&\text{otherwise}.
\end{cases}
\]
For transversal $I$ with $i<t$, the required divisibility is exactly \eqref{eq:admissible} with $m=q-1$, including the total-weight condition at $i=0$. Non-transversal sets have degree zero, and the condition at $i=t$ is automatic, since $\gcd_t(K)=1$.
\end{proof}

\begin{lemma}\label{lem:signed-transversal}
For every sufficiently large $n$ satisfying \eqref{eq:admissible}, there are integers $z_W$, indexed by $W\in\cT_w(n,m)$, such that
\begin{equation}\label{eq:signed-transversal}
\sum_{W\in\cT_w(n,m)}z_W\,\one_{\cT_t(W)}=\one_{\cT_t(n,m)},
\end{equation}
where $\cT_t(W)=\{u\in\cT_t(n,m):u\subseteq W\}$.
\end{lemma}

\begin{proof}
We apply Lemma~\ref{lem:hypergraph-lattice} with $\Psi=\Xi$, $k=w$, $r=t$, $K=\binom{[w]}t$ and $U=G_X$. By Lemma~\ref{lem:transversal-extendable}, $\Xi$ is $(\eta_*,3t^2)$-extendable for $n\geq6wt^2$, where $\eta_*=2^{-3wt^2}$. The complex $\Xi$ has $N_*=mn$ vertices, so for fixed $m,w,t$ and sufficiently large $n$, we have $\eta_*>N_*^{-1/2}$, and $N_*$ exceeds the size threshold in Lemma~\ref{lem:hypergraph-lattice}. By Lemma~\ref{lem:transversal-divisible}, $G_X$ is $K$-divisible. Hence, there are integers $c_f$, indexed by $f\in\Xi_{[w]}$, such that
\[
G_X=\sum_{f\in\Xi_{[w]}}c_f\,\one_{\binom{f([w])}{t}}.
\]
Each $f\in\Xi_{[w]}$ is an injection from $[w]$ into $X$ with transversal image, and its image $f([w])$ is $\widehat W$ for a unique $W\in\cT_w(n,m)$. The vector $\one_{\binom{f([w])}t}$ is the indicator of the $t$-subsets of $\widehat W$, so its entry at $\widehat u$ equals $1$ if $\widehat u\subseteq\widehat W$ and $0$ otherwise. Embeddings with the same image contribute the same vector, and we collect their coefficients into
\[
z_W=\sum_{f\in\Xi_{[w]},\ f([w])=\widehat W}c_f\in\Z
\ \text{ for } W\in\cT_w(n,m).
\]
Now fix $u\in\cT_t(n,m)$. Evaluating the identity above at $\widehat u$ and grouping the terms by their images, we obtain
\[
1=G_X(\widehat u)
=\sum_{W\in\cT_w(n,m),\ \widehat u\subseteq\widehat W}z_W
=\sum_{W\in\cT_w(n,m),\ u\in\cT_t(W)}z_W,
\]
where the last equality uses \eqref{eq:label-transversal}. This is the $u$-coordinate of \eqref{eq:signed-transversal}.
\end{proof}

\begin{lemma}\label{lem:complete-lattice}
If $n$ is sufficiently large and satisfies \eqref{eq:admissible}, then $G\in\Lambda(\Phi^0)$.
\end{lemma}

\begin{proof}
Note that $G(\psi)=G_X(\widehat{\ell(\psi)})=1$ for every $\psi\in\Phi^0_t$. Take the coefficients $z_W$ from Lemma~\ref{lem:signed-transversal}. We lift this signed combination of blocks to molecules of $\Phi^0$. Fix $W=(S,b)$ in $\cT_w(n,m)$, and for each $k\in[m]$ put $j_k=|b^{-1}(k)|$, so that $0\leq j_k\leq w$ and $\sum_kj_k=w$. Since each color class $P_k$ contains $w$ roles, we may choose an active set $J_W\in\binom Rw$ with $|J_W\cap P_k|=j_k$ for every $k$, and a bijection $\iota_W:J_W\to S$ that sends the roles in $J_W\cap P_k$ to the coordinates in $b^{-1}(k)$. Then $b(\iota_W(u))=\kappa(u)$ for $u\in J_W$. For $n\geq L$, we have $|[n]\setminus S|=n-w\geq L-w=|R\setminus J_W|$, so $\iota_W$ extends to an injection $\phi_W:R\to[n]$. Since $\Phi^0$ is complete, $\phi_W\in\Phi^0_R$, and its active word is $x_{J_W,\phi_W}=W$. By \eqref{eq:molecule-word}, we have $v_{J_W,\phi_W}(\psi)=1$ if $\ell(\psi)\subseteq W$ and $v_{J_W,\phi_W}(\psi)=0$ otherwise. Evaluating \eqref{eq:signed-transversal} at $u=\ell(\psi)$ therefore gives
\[
\sum_{W\in\cT_w(n,m)}z_Wv_{J_W,\phi_W}(\psi)
=\sum_{W\in\cT_w(n,m),\ \ell(\psi)\subseteq W}z_W=1=G(\psi)
\]
for every $\psi\in\Phi^0_t$. Hence, $G=\sum_Wz_Wv_{J_W,\phi_W}$ is an integer combination of molecules of $\Phi^0$, that is, $G\in\Lambda(\Phi^0)$.
\end{proof}

To transfer this membership to the random host, we use Keevash's local characterization of the decomposition lattice. Throughout the rest of this subsection, write $Q=\binom Rt$.

\begin{definition}[see {\cite[Definition~5.14]{Keevash2018}}]\label{def:local-lattice}
Let $\Phi$ be a $\Sigma$-adapted $R$-complex. For $U\in\Z^{\Phi_t}$ and a face $\psi'$ of $\Phi$, the \emph{shadow} of $U$ at $\psi'$ is the vector $U^\sharp_{\psi'}\in\Z^Q$ with coordinates
\begin{equation}\label{eq:shadow}
(U^\sharp_{\psi'})_B=\sum_{\psi\in\Phi_B,\ \psi'\subseteq\psi}U_\psi
\ \text{ for } B\in Q.
\end{equation}
For a face $\theta'$ of a template $A_J$, the \emph{template shadow} $\gamma^\sharp_{J,\theta'}\in\Z^Q$ is given by $(\gamma^\sharp_{J,\theta'})_B=\sum_{\theta\in(A_J)_B,\ \theta'\subseteq\theta}\gamma_{J,\theta}$. For a template embedding $(J,\phi)$, write $v^\sharp_{J,\phi}$ for the family of shadows of $v_{J,\phi}$. For an orbit $O\in\Phi/\Sigma$, let $U^\sharp|_O=(U^\sharp_{\psi'})_{\psi'\in O}\in(\Z^Q)^O$, and let
\[
\gamma^\sharp[O]=\Bigl\{v^\sharp_{J,\phi}|_O:(J,\phi)\in\cA(\Phi),\
v^\sharp_{J,\phi}|_O\neq0\Bigr\}.
\]
The \emph{local lattice} is
\[
\cL_\gamma(\Phi)=\Bigl\{U\in\Z^{\Phi_t}:U^\sharp|_O\in
\langle\gamma^\sharp[O]\rangle_{\Z}\text{ for every }O\in\Phi/\Sigma\Bigr\}.
\]
\end{definition}

On each orbit, a single integer combination must match the shadows at all its faces and on all domains $B\in Q$, while the coefficients may vary from one orbit to another. At the empty face, the shadow records the total weight on each domain $B\in Q$. If $|B'|=t$, the only potentially nonzero coordinate of the shadow at a face with domain $B'$ is the $B'$-coordinate, which equals the value of $U$ at that face. If $|B'|>t$, the shadow is zero. For a template embedding $(J,\phi)$ and a face $\psi'$ of $\Phi$, a nonzero summand in the shadow \eqref{eq:shadow} of $v_{J,\phi}$ at $\psi'$ comes from a template face $\theta$ with $\psi'\subseteq\phi\circ\theta$. Restricting $\theta$ to the domain of $\psi'$ gives a template preimage $\theta'$ of $\psi'$, which is unique by the injectivity of $\phi$, and for this preimage the condition $\psi'\subseteq\phi\circ\theta$ is equivalent to $\theta'\subseteq\theta$. Hence, the transport rule \eqref{eq:molecule} gives
\[
(v^\sharp_{J,\phi})_{\psi'}=
\begin{cases}
\gamma^\sharp_{J,\theta'},&\text{if }\psi'=\phi\circ\theta'\text{ for a face }
\theta'\text{ of }A_J,\\
0,&\text{otherwise}.
\end{cases}
\]
In particular, if $U=\sum_{J,\phi}c_{J,\phi}v_{J,\phi}\in\Lambda(\Phi)$, then $U^\sharp|_O=\sum_{J,\phi}c_{J,\phi}v^\sharp_{J,\phi}|_O$ for every orbit $O$, so $\Lambda(\Phi)\subseteq\cL_\gamma(\Phi)$. The following result gives the reverse inclusion.

\begin{lemma}[see {\cite[Lemma~5.19]{Keevash2018}}]\label{lem:local-lattice}
Let $\Phi$ be a $\Sigma$-adapted, $(\eta,3t^2)$-extendable $R$-complex. If $n=|V(\Phi)|$ is sufficiently large and $\eta>n^{-1/2}$, then $\Lambda(\Phi)=\cL_\gamma(\Phi)$.
\end{lemma}

This is \cite[Lemma~5.19]{Keevash2018} with $q=L$, $r=t$ and $D=1$. To compare two hosts with a common skeleton, we describe the generators $\gamma^\sharp[O]$ in terms of the templates. Fix an orbit $O$ of faces of size at most $t$, choose a representative $\psi_0\in O$, and write $B_0=\operatorname{dom}\psi_0$. All faces in $O$ have the same image, so for each $\xi\in O$ the map $r_\xi=\psi_0^{-1}\circ\xi:\operatorname{dom}\xi\to B_0$ is well defined, where $\psi_0^{-1}$ is the inverse on $\im\psi_0$. Equivalently, $r_\xi$ is the unique map with $\xi=\psi_0\circ r_\xi$, and since $\xi$ is obtained from $\psi_0$ by a color-preserving relabelling, $r_\xi$ is a color-preserving bijection. For $J\in\binom Rw$ and $\theta'\in(A_J)_{B_0}$, define $a_{J,\theta'}\in(\Z^Q)^O$ by
\[
(a_{J,\theta'})_\xi=\gamma^\sharp_{J,\theta'\circ r_\xi}
\ \text{ for } \xi\in O,
\]
where $\theta'\circ r_\xi$ is a color-preserving injection with domain $\operatorname{dom}\xi$, hence a face of $A_J$. Let
\[
\mathscr A(O)=\Bigl\{a_{J,\theta'}:J\in\binom Rw,\
\theta'\in(A_J)_{B_0},\ a_{J,\theta'}\neq0\Bigr\}.
\]
For a fixed representative $\psi_0$, this set is determined by the orbit and the templates.

\begin{lemma}\label{lem:local-generators}
Let $\Psi$ be a $\Sigma$-adapted $R$-complex with a nonempty vertex set, and let $O$ be an orbit of its faces of size at most $t$. Then $\gamma^\sharp[O]\subseteq\mathscr A(O)$. If $\Psi$ is $(\eta,1)$-extendable for some $\eta>0$, then $\gamma^\sharp[O]=\mathscr A(O)$.
\end{lemma}

\begin{proof}
Suppose that a full embedding $\phi\in\Psi_R$ maps a template face $\theta'\in(A_J)_{B_0}$ to $\psi_0$, that is, $\phi\circ\theta'=\psi_0$. Then $\xi=\psi_0\circ r_\xi=\phi\circ(\theta'\circ r_\xi)$ for every $\xi\in O$, so by the transport rule for shadows, the shadow of $v_{J,\phi}$ at $\xi$ is $\gamma^\sharp_{J,\theta'\circ r_\xi}$. Hence, $v^\sharp_{J,\phi}|_O=a_{J,\theta'}$.

To prove $\gamma^\sharp[O]\subseteq\mathscr A(O)$, take a nonzero tuple $v^\sharp_{J,\phi}|_O$ and choose $\xi_0\in O$ at which its shadow is nonzero. The transport rule gives a face $\theta_0$ of $A_J$ with $\phi\circ\theta_0=\xi_0$. Then $\theta'=\theta_0\circ r_{\xi_0}^{-1}$ is a color-preserving injection with domain $B_0$, so $\theta'\in(A_J)_{B_0}$, and $\phi\circ\theta'=\xi_0\circ r_{\xi_0}^{-1}=\psi_0$. By the previous paragraph, $v^\sharp_{J,\phi}|_O=a_{J,\theta'}$, which is nonzero and hence belongs to $\mathscr A(O)$.

For the reverse inclusion, assume that $\Psi$ is $(\eta,1)$-extendable and fix $a_{J,\theta'}\in\mathscr A(O)$. Put $C=\theta'(B_0)$. To realize $\phi\circ\theta'=\psi_0$, we must prescribe $\phi|_C=\psi_C$, where $\psi_C=\psi_0\circ(\theta')^{-1}:C\to V(\Psi)$. Since $\theta'$ is a color-preserving injection, it extends to some $\tau\in\Sigma$, and $\psi_C$ is the relabelling of $\psi_0$ by $\tau^{-1}$. Thus, $\psi_C\in\Psi_C$ by $\Sigma$-adaptedness. To extend $\psi_C$ to all of $R$, take $H=R(1)$, whose vertices are the $L$ pairs $(u,1)$ with $u\in R$, and root it at $F=C\times\{1\}$ by setting $f(u,1)=\psi_C(u)$ for $u\in C$. Every face of $H[F]$ is mapped by $f$ to a restriction of $\psi_C$, so $f$ is a $\Psi$-embedding of $H[F]$. Since $\Psi$ is $(\eta,1)$-extendable, there are at least $\eta|V(\Psi)|^{L-|C|}>0$ completions. Choose one completion $g$, and define $\phi(u)=g(u,1)$ for $u\in R$. Applying the embedding condition to the full face $u\mapsto(u,1)$ of $H$ gives $\phi\in\Psi_R$. Since $g$ extends $f$, we have $\phi|_C=\psi_C$, and hence $\phi\circ\theta'=\psi_0$. By the first paragraph, $a_{J,\theta'}=v^\sharp_{J,\phi}|_O\in\gamma^\sharp[O]$.
\end{proof}

The local characterization now allows us to transfer integer feasibility between two hosts with the same faces of size at most $t$.

\begin{lemma}\label{lem:lattice-invariance}
Fix $\eta>0$. For all sufficiently large $n$, the following holds. Let $\Phi$ and $\Phi'$ be $\Sigma$-adapted, $(\eta,3t^2)$-extendable $R$-complexes on $[n]$ such that $\Phi_B=\Phi'_B$ for every $B\subseteq R$ with $|B|\leq t$. Then $\Lambda(\Phi)=\Lambda(\Phi')$ in $\Z^{\Phi_t}=\Z^{\Phi'_t}$.
\end{lemma}

\begin{proof}
We first show that $\cL_\gamma(\Phi)=\cL_\gamma(\Phi')$. The two hosts contain the same face maps of size at most $t$, and $\Sigma$ acts on these maps by the same relabelling of their domains, so their orbits of faces of size at most $t$ coincide. Fix $U\in\Z^{\Phi_t}=\Z^{\Phi'_t}$ and such a common orbit $O$, and choose the same representative $\psi_0$ of $O$ in both hosts. For a face $\xi\in O$ with domain $D$ and a domain $B\in Q$, the $B$-coordinate of the shadow at $\xi$ is, by \eqref{eq:shadow}, the sum of $U_\psi$ over all $\psi\in\Phi_B$ with $D\subseteq B$ and $\psi|_D=\xi$. Since $\Phi'_B=\Phi_B$, the corresponding sums in $\Phi'$ agree term by term. Hence, $U^\sharp|_O$ is the same in both hosts. By Definition~\ref{def:extendability}, both hosts are $(\eta,1)$-extendable, so Lemma~\ref{lem:local-generators} gives $\gamma^\sharp[O]=\mathscr A(O)$ in both hosts. Therefore, the condition $U^\sharp|_O\in\langle\gamma^\sharp[O]\rangle_{\Z}$ holds in $\Phi$ if and only if it holds in $\Phi'$. If $O$ is an orbit of faces of size greater than $t$ in either host, then no domain $B\in Q$ contains the domain of a face in $O$, so every shadow on $O$ is an empty sum, and $U^\sharp|_O=0$ lies in the corresponding lattice. Taking all orbit conditions together, we obtain $U\in\cL_\gamma(\Phi)$ if and only if $U\in\cL_\gamma(\Phi')$.

Since $\eta$ is fixed and positive, for sufficiently large $n$ we have $\eta>n^{-1/2}$, and both hosts meet the size requirement of Lemma~\ref{lem:local-lattice}. Thus, $\Lambda(\Phi)=\cL_\gamma(\Phi)=\cL_\gamma(\Phi')=\Lambda(\Phi')$.
\end{proof}

\begin{proof}[Proof of Lemma~\ref{lem:host-lattice}]
We first show that the complete role complex $\Phi^0$ on $[n]$ is $\Sigma$-adapted and $(2^{-3Lt^2},3t^2)$-extendable for $n\geq6Lt^2$. It contains every injection from a subset of $R$ into $[n]$, and relabelling the domain preserves injectivity, so $\Phi^0$ is $\Sigma$-adapted. Consider a rooted extension $(H,F,f)$ with $H\subseteq R(3t^2)$, and put $v=|V(H)\setminus F|$. Since $R(3t^2)$ has $3t^2$ copies of each of the $L$ roles, we have $|F|+v=|V(H)|\leq3Lt^2$. A completion assigns distinct coordinates from $[n]\setminus f(F)$ to the $v$ remaining vertices, and every such assignment gives a $\Phi^0$-embedding, because its composition with each face of $H$ is an injection. Thus, the number of completions is exactly $\fall{n-|F|}{v}$. For $v\geq1$, the smallest factor in this falling factorial is $n-|F|-v+1\geq n-3Lt^2+1\geq n/2$, so
\[
\fall{n-|F|}{v}\geq(n/2)^v=2^{-v}n^v\geq2^{-3Lt^2}n^v.
\]
When $v=0$, the root map itself is the unique completion, and the same bound holds.

Now let $\eta'=\min\{\eta,2^{-3Lt^2}\}$. For $n\geq6Lt^2$, it follows from the previous paragraph and Definition~\ref{def:extendability} that both $\Phi$ and $\Phi^0$ are $(\eta',3t^2)$-extendable. They are $\Sigma$-adapted and, by \eqref{eq:full-skeleton}, have the same faces of size at most $t$. Thus, for sufficiently large $n$, Lemma~\ref{lem:lattice-invariance} gives $\Lambda(\Phi)=\Lambda(\Phi^0)$ in their common ambient space $\Z^{\Phi_t}=\Z^{\Phi^0_t}$. By Lemma~\ref{lem:complete-lattice}, the divisibility conditions give $G\in\Lambda(\Phi^0)$. Consequently, $G\in\Lambda(\Phi)$.
\end{proof}

\subsection{Proof of Lemma \texorpdfstring{\ref*{lem:regularity}}{3.4}}\label{subsec:fractional}

In this subsection, we prove Lemma~\ref{lem:regularity}. The weights \eqref{eq:weights} are chosen to balance every color composition of a root with $t$ roles, and the full-extension estimate \eqref{eq:full-extension-count} then gives the regularity. The same count of root placements enters both steps, and we compute it once.

\begin{lemma}\label{lem:root-placement}
For $B\in\binom Rt$ and $J\in\binom Rw$, write $b_a=|B\cap P_a|$ and $j_a=|J\cap P_a|$. Let $\Theta(B,J)$ be the set of color-preserving injections $\theta:B\to J$, and put $M_{B,J}=|\Theta(B,J)|$. Then
\[
M_{B,J}=\prod_{a=1}^m\fall{j_a}{b_a},
\]
where a factor is zero if $b_a>j_a$.
\end{lemma}

\begin{proof}
The $b_a$ distinct roles in $B\cap P_a$ must receive distinct images in $J\cap P_a$, which can be done in $\fall{j_a}{b_a}$ ways. The image sets for different colors lie in disjoint color classes, so the counts multiply. In particular, a placement exists exactly when $b_a\leq j_a$ for every $a$.
\end{proof}

\begin{lemma}\label{lem:weight-identity}
Let $B\in\binom Rt$ and put $b_a=|B\cap P_a|$. Then
\begin{equation}\label{eq:identity}
\sum_{J\in\binom Rw}\alpha_J\prod_{a=1}^m\fall{j_a}{b_a}
=\sum_{J\in\binom Rw}\alpha_JM_{B,J}=1.
\end{equation}
\end{lemma}

\begin{proof}
For nonnegative integers $j_1,\ldots,j_m$ with sum $w$, let $\mathbf j=(j_1,\ldots,j_m)$ and $\mathcal J(\mathbf j)=\{J\in\binom Rw:|J\cap P_a|=j_a\text{ for every }a\in[m]\}$. If $j_a<b_a$ for some $a$, then $M_{B,J}=0$ for every $J\in\mathcal J(\mathbf j)$ by Lemma~\ref{lem:root-placement}. Otherwise, put $k_a=j_a-b_a$. There are $\prod_a\binom w{j_a}$ active sets with composition $\mathbf j$, and for $0\leq b_a\leq j_a\leq w$ we have
\[
\frac{\binom w{j_a}\fall{j_a}{b_a}}{\fall w{j_a}}
=\frac1{j_a!}\cdot\frac{j_a!}{(j_a-b_a)!}=\frac1{k_a!}.
\]
Multiplying these factors and using the prefactor in \eqref{eq:alpha}, we obtain
\[
\sum_{J\in\mathcal J(\mathbf j)}\alpha_JM_{B,J}
=\frac1{m^{w-t}}\cdot\frac{(w-t)!}{k_1!\cdots k_m!}.
\]
Since $\sum_ab_a=|B|=t$, summing over all compositions and substituting $k_a=j_a-b_a$ yields
\[
\sum_{J\in\binom Rw}\alpha_JM_{B,J}
=\frac1{m^{w-t}}\sum_{\substack{k_1,\ldots,k_m\geq0\\k_1+\cdots+k_m=w-t}}
\frac{(w-t)!}{k_1!\cdots k_m!}=1,
\]
where the last equality is the multinomial theorem.
\end{proof}

We next estimate the number of embeddings of one template that contribute at a given root. Each placement of the root inside $J$ gives a full-extension problem, and injectivity makes the resulting families of extensions disjoint.

\begin{lemma}\label{lem:active-set-contribution}
Let $\Phi$ be an $R$-complex on $[n]$ satisfying \eqref{eq:full-skeleton} and \eqref{eq:full-extension-count} with error $c$. For $B\in\binom Rt$, $\psi\in\Phi_B$ and $J\in\binom Rw$, put $S_J(\psi)=\sum_{\phi\in\Phi_R}v_{J,\phi}(\psi)$. Then
\[
|S_J(\psi)-D_nM_{B,J}|\leq cD_nM_{B,J}.
\]
\end{lemma}

\begin{proof}
By \eqref{eq:molecule}, a full map $\phi$ contributes to $S_J(\psi)$ precisely when $\phi\circ\theta=\psi$ for some $\theta\in\Theta(B,J)$. Fix such a placement $\theta$ and put $C_\theta=\theta(B)$. The map $\psi_\theta=\psi\circ\theta^{-1}:C_\theta\to[n]$ is injective and, since $|C_\theta|=t$, belongs to $\Phi_{C_\theta}$ by \eqref{eq:full-skeleton}. The equation $\phi\circ\theta=\psi$ is equivalent to $\phi|_{C_\theta}=\psi_\theta$, so $\mathcal E_\theta=\{\phi\in\Phi_R:\phi\circ\theta=\psi\}$ is the family of full extensions of $\psi_\theta$, and $\bigl||\mathcal E_\theta|-D_n\bigr|\leq cD_n$ by \eqref{eq:full-extension-count}. For a given $\phi$, a placement $\theta$ with $\phi\circ\theta=\psi$ is uniquely determined by $\theta=\phi^{-1}\circ\psi$. Consequently, the families $\mathcal E_\theta$ are disjoint, and $S_J(\psi)=\sum_{\theta\in\Theta(B,J)}|\mathcal E_\theta|$. Hence, $(1-c)D_nM_{B,J}\leq S_J(\psi)\leq(1+c)D_nM_{B,J}$.
\end{proof}

\begin{proof}[Proof of Lemma~\ref{lem:regularity}]
Fix $B\in\binom Rt$ and $\psi\in\Phi_B$. For a fixed active set $J$, every full map $\phi\in\Phi_R$ has weight $y_{J,\phi}=\alpha_J/D_n$ by \eqref{eq:weights}, so the total weighted contribution of $J$ at $\psi$ is $(\alpha_J/D_n)S_J(\psi)$. By Lemma~\ref{lem:active-set-contribution},
\[
\Bigl|\frac{\alpha_J}{D_n}S_J(\psi)-\alpha_JM_{B,J}\Bigr|
=\frac{\alpha_J}{D_n}\bigl|S_J(\psi)-D_nM_{B,J}\bigr|
\leq c\,\alpha_JM_{B,J},
\]
which also covers the case $M_{B,J}=0$, where $S_J(\psi)=0$. Using \eqref{eq:identity} twice, we obtain
\[
\Bigl|\sum_{J,\phi}y_{J,\phi}v_{J,\phi}(\psi)-1\Bigr|
=\Bigl|\sum_J\Bigl(\frac{\alpha_J}{D_n}S_J(\psi)-\alpha_JM_{B,J}\Bigr)\Bigr|
\leq c\sum_J\alpha_JM_{B,J}=c,
\]
where the sums are over $J\in\binom Rw$ and $\phi\in\Phi_R$. Thus, the weighted sum at $\psi$ lies between $1-c$ and $1+c$. Since $B$ and $\psi$ were arbitrary, this proves \eqref{eq:regularity}.

It remains to prove the bounds on the weights. Since $L=mw>t$, we have $\fall{n-t}{L-t}\leq n^{L-t}$, so $y_{J,\phi}=\alpha_J/(\rho\fall{n-t}{L-t})\geq a_0n^{t-L}$. For the upper bound, it suffices to show that $\fall{n-t}{L-t}\geq\frac12n^{L-t}$, which holds for sufficiently large $n$, because
\[
\frac{\fall{n-t}{L-t}}{n^{L-t}}=\prod_{j=t}^{L-1}\Bigl(1-\frac jn\Bigr)
\longrightarrow1
\ \text{ as } n\to\infty.
\]
This gives $y_{J,\phi}\leq2\alpha_J/(\rho n^{L-t})\leq a_1n^{t-L}$ uniformly over all template embeddings, which completes the proof.
\end{proof}

\section{Concluding remarks}\label{sec:conclusion}

In this paper, we proved that, for fixed integers $q\geq3$ and $w>t\geq1$, the divisibility conditions \eqref{eq:admissible} characterize the existence of a $GS(t,w,n,q)$ for all sufficiently large $n$, see Corollary~\ref{cor:characterization}. Thus, the Johnson-type bound \eqref{eq:johnson} is attained with equality at every sufficiently large admissible length, which answers Question~\ref{q:LS} of Liu and Shangguan \cite[Question~5.3]{LiuShangguan2025} for all $q\geq3$ and $w>t\geq1$. The systems we obtain have the additional property that the supports of distinct codewords meet in at most $t$ coordinates. We derived Theorem~\ref{thm:labelled} from Keevash's theorem, Lemma~\ref{lem:keevash}, applied to a random selector complex. We conclude the paper with two open questions.

Our proof gives no explicit value of the threshold $n_0(q,w,t)$ in Theorem~\ref{thm:main}, since it relies on the thresholds in Keevash's theorem. Small admissible lengths are also of interest. Recall that the additional support property forces $n\geq t+(q-1)^t(w-t)$. Hence, at smaller admissible lengths, a generalized Steiner system, if it exists, cannot have this property.

\begin{question}\label{q:threshold}
Can one obtain good bounds on the least $n_0$ such that a $GS(t,w,n,q)$ exists for every admissible $n\geq n_0$? Which smaller admissible lengths admit a $GS(t,w,n,q)$, with or without the property that distinct supports meet in at most $t$ coordinates?
\end{question}

Our proof is non-constructive, since it relies on random selectors and on Keevash's theorem. Liu and Shangguan \cite[Question~5.2]{LiuShangguan2025} asked for explicit constructions attaining lower bounds similar to their non-constructive asymptotic bounds. The corresponding question for exact systems is the following.

\begin{question}\label{q:explicit}
Can one give explicit constructions of $GS(t,w,n,q)$ for all sufficiently large admissible $n$?
\end{question}

\paragraph{Declaration of generative AI use.}
The authors used ChatGPT (OpenAI) and Claude (Anthropic) to assist with language editing and the exploration of ideas. The authors reviewed and revised the AI-assisted content and take responsibility for the manuscript and its results.

\end{document}